\documentclass[a4paper,10pt]{article}
\usepackage{stmaryrd}
\usepackage{amsfonts}
\usepackage{bbm}
\usepackage{chemarrow}
\usepackage{amscd}
\usepackage{mathrsfs}
\usepackage{latexsym,amssymb,amsmath,amscd,amscd,amsthm,amsxtra}
\usepackage[dvips]{graphicx}
\usepackage[utf8]{inputenc}
\usepackage[T1]{fontenc}
\usepackage{lmodern}
\usepackage{amssymb}
\usepackage[all]{xy}
\usepackage{nicefrac,mathtools,enumitem}
\usepackage{microtype}
\usepackage{multirow}
\usepackage{xcolor}
\newtheorem{thm}{Theorem}[section]
\newtheorem{lem}[thm]{Lemma}
\newtheorem{cor}[thm]{Corollary}
\newtheorem{pro}[thm]{Proposition}
\newtheorem{ex}[thm]{Example}
\newtheorem{rmk}[thm]{Remark}
\newtheorem{defi}[thm]{Definition}

\newcommand {\emptycomment}[1]{}

\newcommand{\lon }{\,\rightarrow\,}
\newcommand{\be }{\begin{equation}}
\newcommand{\ee }{\end{equation}}

\newcommand{\pf}{\noindent{\bf Proof.}\ }

\newcommand{\g}{\mathfrak g}
\newcommand{\h}{\mathfrak h}

\newcommand{\huaB}{\mathcal{B}}

\newcommand{\huaA}{\mathcal{A}}
\newcommand{\huaL}{\mathcal{L}}
\newcommand{\huaR}{\mathcal{R}}

\newcommand{\huaX}{\mathcal{X}}
\newcommand{\huaY}{\mathcal{Y}}

\newcommand{\huaO}{\mathcal{O}}

\newcommand{\huaZ}{\mathcal{Z}}

\newcommand{\frkg}{\mathfrak g}
\newcommand{\frkh}{\mathfrak h}

\newcommand{\frkp}{\mathfrak p}
\newcommand{\frkq}{\mathfrak q}

\def\qed{\hfill ~\vrule height6pt width6pt depth0pt}

\newcommand{\re}{\mathrm{Re}}
\newcommand{\aff}{\mathrm{aff}}

\newcommand{\Id}{\rm{Id}}

\newcommand{\br}[1]{   [ \cdot,    \cdot  ]   }

\newcommand{\Hom}{\mathrm{Hom}}

\newcommand{\gl}{\mathfrak {gl}}

\newcommand{\ad}{\mathrm{ad}}

\begin{document}
\title{
{Symplectic, product and complex structures on 3-Lie algebras
} }
\author{Yunhe Sheng and  Rong Tang
 }
\date{}
\maketitle

\begin{abstract}
In this paper, first we introduce the notion of a phase space of a 3-Lie algebra and show that a 3-Lie algebra has a phase space if and only if it is sub-adjacent to a 3-pre-Lie algebra. Then we introduce the notion of a product structure on a 3-Lie algebra using the Nijenhuis condition as the integrability condition. A 3-Lie algebra enjoys a product structure if and only if it is the direct sum (as vector spaces) of two subalgebras. We find that there are four types special integrability conditions, and each of them gives rise to a special decomposition of the original 3-Lie algebra. They are also related to $\huaO$-operators, Rota-Baxter operators and matched pairs of 3-Lie algebras. Parallelly, we introduce the notion of a complex structure on a 3-Lie algebra and there are also four types special integrability conditions. Finally, we add compatibility conditions between a complex structure and a product structure, between a symplectic structure and a paracomplex structure, between a symplectic structure and a complex structure, to introduce the notions of a complex product structure, a para-K\"{a}hler structure and a pseudo-K\"{a}hler structure on a 3-Lie algebra. We use 3-pre-Lie algebras to construct these structures. Furthermore,   a Levi-Civita product is introduced associated to a pseudo-Riemannian 3-Lie algebra and deeply studied.
\end{abstract}

\tableofcontents

\setcounter{section}{0}

\section{Introduction}

A symplectic structure on a Lie algebra $\g$ is a nondegenerate 2-cocycle $\omega\in\wedge^2\g^*$. The underlying structure of a symplectic Lie algebra is a quadratic pre-Lie algebra \cite{Chu}. An almost product structure on a Lie algebra $\g$ is a linear map $E$ satisfying $E^2=\Id$. If in addition, $E$ also satisfies the following integrability condition
$$
[Ex,Ey]=E([Ex,y]+[x,Ey]-E[x,y]),\quad \forall x,y\in\g,
$$
then $E$ is called a product structure. The above integrability condition is called the Nijenhuis condition. An equivalent characterization of a product structure is that $\g$ is the direct sum (as vector spaces) of two subalgebras. An almost complex structure on a Lie algebra $\g$ is a linear map $J$ satisfying $J^2=-\Id$. A  complex structure on a Lie algebra is an almost complex structure that satisfies the Nijenhuis condition. Adding compatibility conditions between a complex structure and a product structure, between a symplectic structure and a paracomplex structure, between a symplectic structure and a complex structure, one obtains a complex product structure, a paraK\"{a}hler structure and a pseudo-K\"{a}hler structure respectively.  These structures play important roles in  algebra, geometry and mathematical physics, and are widely studied. See \cite{Alek,Andrada0,ABD,ABDO,AS,Baibialgebra,Bai-2,Banayadi0,Benayadi,Calvaruso0,Calvaruso,Poon1,Poon2,Li,Salamon} for more details.

Generalizations of Lie algebras to higher arities, including 3-Lie
algebras and more generally, $n$-Lie
algebras~\cite{Filippov,Kasymov,Tcohomology}, have attracted attention
from several fields of mathematics and physics. It is the algebraic
structure corresponding to Nambu mechanics \cite{Gautheron,N,T}. In
particular, the study of 3-Lie algebras plays an important role in
string theory. In \cite{Basu},
 Basu and Harvey suggested to replace the Lie algebra appearing in the Nahm equation by a 3-Lie algebra for the  lifted Nahm equations.
Furthermore, in the context of Bagger-Lambert-Gustavsson model of multiple
M2-branes, Bagger-Lambert  managed to construct, using a ternary bracket,  an $N=2$
 supersymmetric version of the worldvolume theory of the M-theory membrane, see \cite{BL0}. An extensive literatures are related to this pioneering work, see \cite{BL3,BL2,HHM,P}. See the review article \cite{review} for more details. In particular,   metric 3-algebras were deeply studied in the seminal works \cite{DFM,DFMR,DFMR2}. In \cite{Liu-Sheng-Bai-Chen}, the authors introduced the notion a Nijenhuis operator on an $n$-Lie algebra, which generates a trivial deformation.

 The purpose of this paper is to   study   symplectic structures, product structure and complex structures on 3-Lie algebras and these combined structures. In the case of Lie algebras, pre-Lie algebras play important roles in these studies. It is believable that 3-pre-Lie algebras will play important roles in the corresponding studies. Thus, first we introduce the notion of a representation of a 3-pre-Lie algebra and construct the associated semidirect product 3-pre-Lie algebra. Several important properties of representations of 3-pre-Lie algebras are studied.

 Note that the notion of a symplectic structure on a 3-Lie algebra was introduced in \cite{BGS} and it is shown that the underlying structure of a symplectic 3-Lie algebra is a quadratic 3-pre-Lie algebra. We introduce the notion of a phase space of a 3-Lie algebra $\g$, which is a symplectic 3-Lie algebra $\g\oplus \g^*$ satisfying some conditions, and show that a 3-Lie algebra has a phase space if and only if it is sub-adjacent to a 3-pre-Lie algebra. We also introduce the notion of a Manin triple of 3-pre-Lie algebras and show that there is a one-to-one correspondence between Manin triples of 3-pre-Lie algebras and phase spaces of   3-Lie algebras.

 An almost product structure on a 3-Lie algebra $\g$ is defined to be a linear map $E:\g\longrightarrow\g$ satisfying $E^2=\Id$. It is challengeable to add an integrability condition on an almost product structure to obtain a product structure on a 3-Lie algebra. We note that the Nijenhuis condition (see \eqref{eq:Nejenhuiscon}) given in \cite{Liu-Sheng-Bai-Chen} is the correct integrability condition. Let us explain this issue. Denote by $\g_{\pm}$ the eigenspaces corresponding to eigenvalues $\pm1$ of an almost product structure $E$. Then it is obvious that $\g=\g_+\oplus \g_-$ as vector spaces. The Nijenhuis condition ensures that both $\g_+$ and $\g_-$ are subalgebras. This is what ``integrability'' means. Moreover, we find that there are four types special integrability conditions, which are called strict product structure, abelian product structure, strong abelian product structure and perfect product structure respectively, each of them gives rise to a special decomposition of the original 3-Lie algebra. See the following table for a precise description:
\begin{tabular}{|c|c|c|}
\hline product   & $E[x,y,z]_\g=[Ex,y,z]_\g+[x,Ey,z]_\g+[x,y,Ez]_\g$ & $\g=\g_+\oplus\g_-$ \\
structure&$-E([Ex,Ey,z]_\g+[Ex,y,Ez]_\g+[x,Ey,Ez]_\g)$&$[\g_+,\g_+,\g_+]_\g\subset\g_+$\\
&$+[Ex,Ey,Ez]_\g$&$[\g_-,\g_-,\g_-]_\g\subset\g_-$\\\hline
 strict product& $E[x,y,z]_\g=[Ex,y,z]_\g$ & $[\g_+,\g_+,\g_-]_\g=0$\\
structure & {} & $[\g_-,\g_-,\g_+]_\g=0$\\\hline
 abelian product  & $[x,y,z]_\g=-[x,Ey,Ez]_\g-[Ex,y,Ez]_\g-[Ex,Ey,z]_\g$ & $[\g_+,\g_+,\g_+]_\g=0$\\
 structure  & {} & $[\g_-,\g_-,\g_-]_\g=0$\\\hline
    & $[x,y,z]_\g=E[Ex,y,z]_\g+E[x,Ey,z]_\g+E[x,y,Ez]_\g$ & $[\g_+,\g_+,\g_+]_\g=0$\\
strong abelian & {} & $[\g_-,\g_-,\g_-]_\g=0$\\
 product structure  & {$\huaO$-operators} & $[\g_+,\g_+,\g_-]_\g\subset\g_+$\\
{} & {Rota-Baxter operators} & $[\g_-,\g_-,\g_+]_\g\subset\g_-$\\\hline
 perfect product  & $E[x,y,z]_\g=[Ex,Ey,Ez]_\g$ & $[\g_+,\g_+,\g_-]_\g\subset\g_-$\\
{structure} & {involutive automorphisms} & $[\g_-,\g_-,\g_+]_\g\subset\g_+$\\\hline
\end{tabular}
It is surprised that a strong abelian product structure is also an $\huaO$-operator   on a 3-Lie algebra associated to the adjoint representation. Since an $\huaO$-operator   on a 3-Lie algebra associated to the adjoint representation is also a Rota-Baxter operator \cite{RB3Lie,PBG}, it turns out that involutive Rota-Baxter operator can also serve as an integrability condition. This is totally different from the case of Lie algebras. Furthermore, by the definition of a perfect product structure, an involutive automorphism of a 3-Lie algebra can also serve as an integrability condition. This is also a new phenomenon. Note that the decomposition that a perfect product structure gives is exactly the condition required in the definition of a matched pair of 3-Lie algebras \cite{BGS}. Thus, this kind of product structure will be frequently used in our studies.

An   almost complex structure on a 3-Lie algebra $\g$ is defined to be a linear map $J:\g\longrightarrow\g$ satisfying $J^2=-\Id$. With the above motivation, we define a complex structure on a 3-Lie algebra $\g$ to be an   almost complex structure satisfying the Nijenhuis condition. Then $\g_i$ and $\g_{-i}$, which are eigenspaces of eigenvalues $\pm i$ of a complex linear map $J_{\mathbb C}$ (the complexification of $J$) are subalgebras of the 3-Lie algebra $\g_{\mathbb C}$, the complexification of $\g$. Parallel to the case of product structures, there are also four types special integrability conditions, and each of them gives rise to a special decomposition of $\g_{\mathbb C}$:
\begin{tabular}{|c|c|c|}
\hline complex   & $J[x,y,z]_\g=[Jx,y,z]_\g+[x,Jy,z]_\g+[x,y,Jz]_\g$ & $\g_{\mathbb C}=\g_i\oplus\g_{-i}$ \\
structure&$+J([Jx,Jy,z]_\g+[Jx,y,Jz]_\g+[x,Jy,Jz]_\g)$&$[\g_i,\g_i,\g_i]_{\g_{\mathbb C}}\subset\g_i$\\
&$-[Jx,Jy,Jz]_\g$&$[\g_{-i},\g_{-i},\g_{-i}]_{\g_{\mathbb C}}\subset\g_{-i}$\\\hline
 strict complex& $J[x,y,z]_\g=[Jx,y,z]_\g$ & $[\g_i,\g_i,\g_{-i}]_{\g_{\mathbb C}}=0$\\
structure & {} & $[\g_{-i},\g_{-i},\g_i]_{\g_{\mathbb C}}=0$\\\hline
 abelian complex  & $[x,y,z]_\g=[x,Jy,Jz]_\g+[Jx,y,Jz]_\g+[Jx,Jy,z]_\g$ & $[\g_i,\g_i,\g_i]_{\g_{\mathbb C}}=0$\\
 structure  & {} & $[\g_{-i},\g_{-i},\g_{-i}]_{\g_{\mathbb C}}=0$\\\hline
    & $[x,y,z]_\g=-J([Jx,y,z]_\g+[x,Jy,z]_\g+[x,y,Jz]_\g)$ & $[\g_i,\g_i,\g_i]_{\g_{\mathbb C}}=0$\\
strong abelian & {} & $[\g_{-i},\g_{-i},\g_{-i}]_{\g_{\mathbb C}}=0$\\
complex structure  & {$\huaO$-operators} & $[\g_i,\g_i,\g_{-i}]_{\g_{\mathbb C}}\subset\g_i$\\
{} & {Rota-Baxter operators} & $[\g_{-i},\g_{-i},\g_i]_{\g_{\mathbb C}}\subset\g_{-i}$\\\hline
 perfect complex  & $J[x,y,z]_\g=-[Jx,Jy,Jz]_\g$ & $[\g_i,\g_i,\g_{-i}]_{\g_{\mathbb C}}\subset\g_{-i}$\\
{structure} & {anti-involutive automorphisms} & $[\g_{-i},\g_{-i},\g_i]_{\g_{\mathbb C}}\subset\g_i$\\\hline
\end{tabular}

Then we   add a compatibility condition between a complex structure and a product structure on a 3-Lie algebra to define a complex product structure on a 3-Lie algebra. We give an equivalent characterization of a complex product structure on a 3-Lie algebra $\g$ using the decomposition of $\g$. We add   a compatibility condition between a symplectic structure and a paracomplex structure on a 3-Lie algebra to define a paraK\"{a}hler structure on a 3-Lie algebra. An equivalent characterization of a paraK\"{a}hler structure on a 3-Lie algebra $\g$ is also given using the decomposition of $\g$. Associated to a paraK\"{a}hler structure on a 3-Lie algebra, there is also a pseudo-Riemannian structure. We introduce the notion of a Livi-Civita product associated to a pseudo-Riemannian 3-Lie algebra, and give its precise formulas. Finally, we add a compatibility condition between a symplectic structure and a complex   structure on a 3-Lie algebra to define a pseudo-K\"{a}hler structure on a 3-Lie algebra. The relation between a paraK\"{a}hler structure and a pseudo-K\"{a}hler structure on a 3-Lie algebra is investigated.  We construct complex product structures, paraK\"{a}hler structures and   pseudo-K\"{a}hler structures in terms of 3-pre-Lie algebras. We also give examples of symplectic structures, product structures, complex structures, complex product structure, paraK\"{a}hler structures and   pseudo-K\"{a}hler structures on the $4$-dimensional Euclidean $3$-Lie algebra $A_{4}$ given in \cite{BL0}.

The paper is organized as follows. In Section 2, we recall Nijenhuis operators on 3-Lie algebras and 3-pre-Lie algebras. In Section 3, we study representations of 3-pre-Lie algebras. In Section 4, we introduce the notion of a phase space of a 3-Lie algebra and show that a 3-Lie algebra has a phase space if and only if it is sub-adjacent to a 3-pre-Lie algebra. We also introduce the notion of a Manin triple of 3-pre-Lie algebras and study its relation with phase spaces of 3-Lie algebras. In Section 5, we introduce the notion of a product structure on a 3-Lie algebra and give four special integrability conditions. In Section 6, we introduce the notion of a complex structure on a 3-Lie algebra and give four special integrability conditions. In Section 7, we introduce the notion of a complex product structure on a 3-Lie algebra and give its equivalent characterization. In Section 8, we introduce the notion of a paraK\"{a}hler structure on a 3-Lie algebra and give its equivalent characterization. Moreover, we give a detailed study on the associated Levi-Civita product. In Section 9, we  introduce the notion of a pseudo-K\"{a}hler structure on a 3-Lie algebra and study the relation with a paraK\"{a}hler structure.

In this paper, we work over the real field $\mathbb R$ and the complex field $\mathbb C$ and all the vector
spaces are finite-dimensional.

\vspace{2mm}
 \noindent {\bf Acknowledgement:} We give our warmest thanks to Chengming Bai for very useful comments and discussions.  This research is supported by NSFC (11471139) and NSF of Jilin Province (20170101050JC).

  \section{Preliminaries}
In this section, first we recall the notion of a Nijenhuis operator on a 3-Lie algebra, which will be frequently used as the integrability condition in our later studies. Then  we recall the notion of a  3-pre-Lie algebra, which is the main tool to construct examples of symplectic, product and complex structures on 3-Lie algebras.
\begin{defi}\label{defi of n-LA}
A {\bf $3$-Lie algebra} is a vector space $\g$ together with a trilinear  skew-symmetric bracket $[\cdot,\cdot,\cdot]_\g:\wedge^3\g\longrightarrow\g$ such that   the following {\bf fundamental} identity holds:
\begin{eqnarray}\label{FI}
[x,y,[z,w,v]_\g]_\g=[[x,y,z]_\g,w,v]_\g+[z,[x,y,w]_\g,v]_\g+[z,w,[x,y,v]_\g]_\g,\quad \forall x,y,z,w,v\in\g.
\end{eqnarray}
\end{defi}
\emptycomment{
\begin{defi}
A  morphism of $3$-Lie algebras $f:(\g,[\cdot,\cdot,\cdot]_\g)\lon(\h,[\cdot,\cdot,\cdot]_\h)$ is a linear map $f:\g\lon\h$ such that
\begin{eqnarray}
f[x,y,z]_\g=[f(x),f(y),f(z)]_\h,\hspace{3mm}\forall x,y,z\in \mathfrak{g}.
\end{eqnarray}
\end{defi}
}
For $x,y\in\g$, define $\ad:\wedge^{2}\g\longrightarrow\gl(\g)$ by
$$\ad_{x,y}z=[x,y,z]_\g, \quad\forall z\in \g.$$
Then \eqref{FI} is equivalent to that $\ad_{x,y}$ is a derivation, i.e.
$$
\ad_{x,y}[z,w,v]_\g=[\ad_{x,y}z,w,v]_\g+[z,\ad_{x,y}w,v]_\g+[z,w,\ad_{x,y}v]_\g,\quad \forall x,y\in\g.
$$

Let $(\g,[\cdot,\cdot,\cdot]_\g)$ be a $3$-Lie algebra, and $N:\g\longrightarrow\g$ a linear map. Define a $3$-ary bracket $[\cdot,\cdot,\cdot]_N^1:\wedge^3\g\longrightarrow\g$  by
   \begin{equation}\label{eq:bracket(1)}
    [x,y,z]_N^{1}=[Nx,y,z]_\g+[x,Ny,z]_\g+[x,y,Nz]_\g-N[x,y,z]_\g.
  \end{equation}
  Then we define $3$-ary bracket $[\cdot,\cdot,\cdot]_N^2:\wedge^3\g\longrightarrow\g$ by
    \begin{equation}\label{eq:bracket (j)}
    [x,y,z]_N^{2}=[Nx,Ny,z]_\g+[x,Ny,Nz]_\g+[Nx,y,Nz]_\g-N[x,y,z]_N^{1}.
  \end{equation}

\begin{defi}\label{defi:Nijenhuis}{\rm (\cite{Liu-Sheng-Bai-Chen})}
Let $(\g,[\cdot,\cdot,\cdot]_\g)$ be a $3$-Lie algebra. A linear map    $N:\g\longrightarrow\g$   is called a {\bf Nijenhuis operator} if the following {\bf Nijenhuis condition} is satisfied:
\begin{equation}\label{eq:Nijenhuis(n)}
[Nx,Ny,Nz]_\g=N[x,y,z]_N^{2},\quad\forall x,y,z\in \g.
\end{equation}
\end{defi}

More precisely,   a linear map $N:\g\longrightarrow\g$ of a $3$-Lie algebra $(\g,[\cdot,\cdot,\cdot]_\g)$ is a Nijenhuis operator if and only if
\begin{eqnarray}
\nonumber[Nx,Ny,Nz]_\g&=&N[Nx,Ny,z]_\g+N[x,Ny,Nz]_\g+N[Nx,y,Nz]_\g\\\nonumber
&&-N^2[Nx,y,z]_\g-N^2[x,Ny,z]_\g-N^2[x,y,Nz]_\g\\
\label{eq:Nejenhuiscon}&&+N^3[x,y,z]_\g.
\end{eqnarray}

\begin{defi}{\rm (\cite{Kasymov})}\label{defi:usualrep}
 A {\bf representation}   of a $3$-Lie algebra $(\g,[\cdot,\cdot,\cdot]_\g)$ on a vector space  $V$ is a linear map $\rho:\wedge^2\frkg\longrightarrow \gl( V),$ such that for all $x_1,x_2,x_3,x_4\in\g,$ there holds:
\begin{eqnarray*}
 &\rho([x_1,x_2,x_3]_\g, x_4) +\rho(x_3,[x_1,x_2, x_4]_\g) =[\rho(x_1,x_2),\rho(x_3,x_4)];\\
&\rho([x_1,x_2,x_3]_\g,x_4)=\rho(x_1,x_2)\circ\rho(x_3,x_4)+\rho(x_2,x_3)\circ\rho(x_1,x_4)+\rho(x_3,x_1)\circ\rho(x_2,x_4).
\end{eqnarray*}
\end{defi}

\begin{ex}{\rm
Let $(\g,[\cdot,\cdot,\cdot]_\g)$ be a $3$-Lie algebra. The linear map $\ad:\wedge^{2}\g\longrightarrow\gl(\g)$ defines a representation of the $3$-Lie algebra $\g$ on itself, which we call the {\bf adjoint representation} of $\g$.
}
\end{ex}
Let $A$ be a vector space. For a linear map $\phi:A\otimes A\lon\gl(V)$, we define a linear map $\phi^*: A\otimes A\lon\gl(V^*)$ by
\begin{eqnarray*}
\langle \phi^*(x,y)\alpha,v\rangle=-\langle\alpha, \phi(x,y)v\rangle,\,\,\,\,\forall \alpha\in V^*,x,y\in\g,v\in V.
\end{eqnarray*}

\emptycomment{
Let $(V,\rho)$ be a representation of a $3$-Lie algebra $(\g,[\cdot,\cdot,\cdot]_\g)$. Define $\rho^*:\wedge^2\frkg\longrightarrow \gl(V^*)$ by
\begin{eqnarray}
\langle\rho^*(x_1,x_2)\alpha,v\rangle=-\langle\alpha,\rho(x_1,x_2)v\rangle,\,\,\,\,\forall \alpha\in V^*,x_1,x_2\in\g,v\in V.
\end{eqnarray}
}

\begin{lem}{\rm (\cite{BGS})}\label{dual-rep-3-Lie}
Let $(V,\rho)$ be a representation of a $3$-Lie algebra $(\g,[\cdot,\cdot,\cdot]_\g)$. Then $(V^*,\rho^*)$ is a representation of the $3$-Lie algebra $(\g,[\cdot,\cdot,\cdot]_\g)$, which is called the {\bf dual representation}.
\end{lem}

\begin{lem}\label{lem:semidirectp}
Let $\g$ be a $3$-Lie algebra, $V$  a vector space and $\rho:
\wedge^2\g\rightarrow \gl(V)$  a skew-symmetric linear
map. Then $(V;\rho)$ is a representation of $\g$ if and only if there
is a $3$-Lie algebra structure (called the {\bf semidirect product})
on the direct sum of vector spaces  $\g\oplus V$, defined by
\begin{equation}\label{eq:sum}
[x_1+v_1,x_2+v_2,x_3+v_3]_{\rho}=[x_1,x_2,x_3]_\g+\rho(x_1,x_2)v_3+\rho(x_2,x_3)v_1+\rho(x_3,x_1)v_2,
\end{equation}
for all $x_i\in \g, v_i\in V, 1\leq i\leq 3$. We denote this semidirect product $3$-Lie algebra by $\g\ltimes_\rho V.$
\end{lem}

\begin{defi}
Let $A$ be a vector space with a linear map $\{\cdot,\cdot,\cdot\}:\otimes^3 A\lon A$. The pair $(A,\{\cdot,\cdot,\cdot\})$ is called a $3$-{\bf pre-Lie algebra} if the following identities hold:
\begin{eqnarray}
\{x,y,z\}                                              &=&-\{y,x,z\}\\
\nonumber\{x_1,x_2,\{x_3,x_4,x_5\}\}                   &=&\{[x_1,x_2,x_3]_C,x_4,x_5\}+\{x_3,[x_1,x_2,x_4]_C,x_5\}\\
                                                       &&+\{x_3,x_4,\{x_1,x_2,x_5\}\}\\
\nonumber \{[x_1,x_2,x_3]_C,x_4,x_5\}                  &=&\{x_1,x_2,\{x_3,x_4,x_5\}\}+\{x_2,x_3,\{x_1,x_4,x_5\}\}\\
                                                       &&+\{x_3,x_1,\{x_2,x_4,x_5\}\},
\end{eqnarray}
where $x,y,z,x_i\in A,1\le i\le 5$ and $[\cdot,\cdot,\cdot]_C$ is defined by
\begin{eqnarray}
[x,y,z]_C\triangleq \{x,y,z\}+\{y,z,x\}+\{z,x,y\},\,\,\,\,\forall x,y,z\in A.
\end{eqnarray}
\end{defi}

\begin{pro}{\rm (\cite[Proposition 3.21]{BGS})}\label{3-pre-Lie}
 Let $(A,\{\cdot,\cdot,\cdot\})$ be a $3$-pre-Lie algebra. Then $(A,[\cdot,\cdot,\cdot]_C)$ is a $3$-Lie algebra, which is called the sub-adjacent $3$-Lie algebra of $A$, and denoted by $A^c$.  $(A,\{\cdot,\cdot,\cdot\})$ is called the compatible $3$-pre-Lie algebra structure on the $3$-Lie algebra $A^c$.
 \end{pro}

Define the left multiplication $L:\wedge^2 A\longrightarrow\gl(A)$ by $L(x,y)z=\{x,y,z\}$ for all $x,y,z\in A$. Then $(A,L)$ is a representation of the $3$-Lie algebra $A^c$. Moreover, we define the right multiplication $R:\otimes^2 A\lon\gl(A)$ by $R(x,y)z=\{z,x,y\}$. If there is a $3$-pre-Lie algebra structure on its dual space $A^*$, we denote the left multiplication and right multiplication by $\huaL$ and $\huaR$ respectively.

\emptycomment{
\begin{defi}
A  morphism of $3$-pre-Lie algebras $f:(A,\{\cdot,\cdot,\cdot\}_A)\lon(A',\{\cdot,\cdot,\cdot\}_{A'})$ is a linear map $f:A\lon A'$ such that
\begin{eqnarray}
f\{x,y,z\}_A)=\{f(x),f(y),f(z)\}_{A'}\hspace{3mm}\forall x,y,z\in A.
\end{eqnarray}
\end{defi}
}

\begin{defi}{\rm (\cite[Definition 3.16]{BGS})}\label{3-Lie-O-operator}
Let $(\g,[\cdot,\cdot,\cdot]_\g)$ be a $3$-Lie algebra and $(V,\rho)$ a representation. A linear operator $T:V\lon\g$ is called an {\bf$\huaO$-operator} associated to $(V,\rho)$ if $T$ satisfies:
\begin{eqnarray}
[Tu,Tv,Tw]_\g=T(\rho(Tu,Tv)w+\rho(Tv,Tw)u+\rho(Tw,Tu)v),\,\,\,\,\forall u,v,w\in V.
\end{eqnarray}
\end{defi}

\emptycomment{
\begin{pro}
Let $(\g,[\cdot,\cdot,\cdot]_\g)$ be a $3$-Lie algebra and $(V,\rho)$ a representation. Suppose that the linear map $T:V\lon\g$ is an $\huaO$-operator associated to $(V,\rho)$. Then there exists a $3$-pre-Lie algebra structure on $V$ given by
\begin{eqnarray}
\{u,v,w\}=\rho(Tu,Tv)w,\,\,\,\,\forall u,v,w\in V.
\end{eqnarray}
\end{pro}
}

\begin{pro}{\rm (\cite[Proposition 3.27]{BGS})}\label{3-Lie-compatible-3-pre-Lie}
Let $(\g,[\cdot,\cdot,\cdot]_\g)$ be a $3$-Lie algebra. Then there is a compatible $3$-pre-Lie algebra if and only if there exists an invertible $\huaO$-operator $T:V\lon\g$ associated to a representation $(V,\rho)$. Furthermore, the compatible $3$-pre-Lie structure on $\g$ is given by
\begin{eqnarray}
\{x,y,z\}=T\rho(x,y)T^{-1}(z),\,\,\,\,\forall x,y,z\in \g.
\end{eqnarray}
\end{pro}

\emptycomment{
\begin{lem}
Given a matched pair $(\g,\h)$ of $3$-Lie algebras, there is a $3$-Lie algebra structure $\g\bowtie\h$ on the direct sum vector space $\g\oplus\h$ with the bracket
\begin{eqnarray}
\nonumber[x_1+u_1,x_2+u_2,x_3+u_3]_{\g\bowtie\h}&=&[x_1,x_2,x_3]_{\g}+\nu_\h(u_1)(x_2,x_3)+\nu_\h(u_2)(x_3,x_1)+\nu_\h(u_3)(x_1,x_2)\\
                                       &&+\nonumber\rho_\h(u_1,u_2)x_3+\rho_\h(u_2,u_3)x_1+\rho_\h(u_3,u_1)x_2\\
                                       &&+\nonumber[u_1,u_2,u_3]_\h+\nu_\g(x_1)(u_2,u_3)+\nu_\g(x_2)(u_3,u_1)+\nu_\g(x_3)(u_1,u_2)\\                      &&+\rho_\g(x_1,x_2)u_3+\rho_\g(x_2,x_3)u_1+\rho_\g(x_3,x_1)u_2.
\end{eqnarray}
Conversely, if $\g\oplus\h$ has a $3$-Lie algebra structure for which $\g$ and $\h$ are $3$-Lie algebras, then the four linear maps defined by
\begin{eqnarray}
&&\rho_{\g}(x,y)u=P_{\h}([x,y,u]_{\g\oplus\h}),\,\,\,\,\nu_\g(x)(u,v)=P_{\h}([x,u,v]_{\g\oplus\h}),\\
&&\rho_{\h}(u,v)x=P_{\g}([u,v,x]_{\g\oplus\h}),\,\,\,\,\nu_\h(u)(x,y)=P_{\g}([u,x,y]_{\g\oplus\h}),
\end{eqnarray}
where $P_{\g}$ and $P_{\h}$ are the natural projection of $\g\oplus\h$ to $\g$ and $\h$ respectively. Moreover, they endow the couple $(\g,\h)$ with a structure of a matched pair.  Thus, the three Lie algebras
$(\g\oplus\h, \g, \h)$ form a double $3$-Lie algebra.
\end{lem}
}

\section{Representations of 3-pre-Lie algebras}

In this section, we introduce the notion of a representation of a 3-pre-Lie algebra, construct the corresponding semidirect product 3-pre-Lie algebra and give the dual representation.
\begin{defi}\label{defi:rep3-pre-Lie}
  A {\bf representation} of a $3$-pre-Lie algebra $(A,\{\cdot,\cdot,\cdot\})$   on a vector space $V$ consists of a pair $(\rho,\mu)$, where $\rho:\wedge^2 A\rightarrow \gl(V)$ is a representation of the $3$-Lie algebra $A^c$ on $V$ and $\mu:\otimes^2 A\rightarrow \gl(V)$ is a linear map such that  for all $x_1,x_2,x_3,x_4\in A$, the following equalities hold:
\begin{eqnarray}
\nonumber \rho(x_1,x_2)\mu(x_3,x_4) &=&\mu(x_3,x_4)\rho(x_1,x_2)-\mu(x_3,x_4)\mu(x_2,x_1)\\
                          \label{rep1}&&+\mu(x_3,x_4)\mu(x_1,x_2)+\mu([x_1,x_2,x_3]_C,x_4)+\mu(x_3,\{x_1,x_2,x_4\}),\\
\label{rep2} \mu([x_1,x_2,x_3]_C,x_4)&=&\rho(x_1,x_2)\mu(x_3,x_4)+\rho(x_2,x_3)\mu(x_1,x_4)+\rho(x_3,x_1)\mu(x_2,x_4),\\
\nonumber \mu(x_1,\{x_2,x_3,x_4\}) &=&\mu(x_3,x_4)\mu(x_1,x_2)+\mu(x_3,x_4)\rho(x_1,x_2)\\
                         \nonumber &&-\mu(x_3,x_4)\mu(x_2,x_1)-\mu(x_2,x_4)\mu(x_1,x_3)\\
                         \label{rep3} &&-\mu(x_2,x_4)\rho(x_1,x_3)+\mu(x_2,x_4)\mu(x_3,x_1)+\rho(x_2,x_3)\mu(x_1,x_4),\\
\nonumber \mu(x_3,x_4)\rho(x_1,x_2) &=&\mu(x_3,x_4)\mu(x_2,x_1)-\mu(x_3,x_4)\mu(x_1,x_2)\\
                        \label{rep4}  &&+\rho(x_1,x_2)\mu(x_3,x_4)-\mu(x_2,\{x_1,x_3,x_4\})+\mu(x_1,\{x_2,x_3,x_4\}).
\end{eqnarray}
\end{defi}

\emptycomment{
\begin{eqnarray}
\nonumber L(x_1,x_2)(R(x_4,x_5)v_3) &=&R(x_4,x_5)(L(x_1,x_2)v_3)-R(x_4,x_5)(R(x_2,x_1)v_3)\\
                         \nonumber &&+R(x_4,x_5)(R(x_1,x_2)v_3)+R([x_1,x_2,x_4]_C,x_5)v_3\\
                          \label{rep1}&&+R(x_4,\{x_1,x_2,x_5\})v_3,\\
\nonumber -R([x_1,x_2,x_3]_C,x_5)v_4&=&-L(x_1,x_2)(R(x_3,x_5)v_4)-L(x_2,x_3)(R(x_1,x_5)v_4)\\
                         \label{rep2}&&-L(x_3,x_1)(R(x_2,x_5)v_4),\\
\nonumber R(x_2,\{x_3,x_4,x_5\})v_1 &=&R(x_4,x_5)(R(x_2,x_3)v_1)+R(x_4,x_5)(L(x_2,x_3)v_1)\\
                         \nonumber &&-R(x_4,x_5)(R(x_3,x_2)v_1)-R(x_3,x_5)(R(x_2,x_4)v_1)\\
                         \nonumber &&-R(x_3,x_5)(L(x_2,x_4)v_1)+R(x_3,x_5)(R(x_4,x_2)v_1)\\
                         \label{rep3}&&+L(x_3,x_4)(R(x_2,x_5)v_1),\\
\nonumber R(x_4,x_5)(L(x_1,x_2)v_3) &=&R(x_4,x_5)(R(x_2,x_1)v_3)-R(x_4,x_5)(R(x_1,x_2)v_3)\\
                         \nonumber &&+L(x_1,x_2)(R(x_4,x_5)v_3)-R(x_2,\{x_1,x_4,x_5\})v_3\\
                         \label{rep4} &&+R(x_1,\{x_2,x_4,x_5\})v_3.
\end{eqnarray}
}

Let $(A,\{\cdot,\cdot,\cdot\})$ be a $3$-pre-Lie algebra and $\rho$   a representation of the sub-adjacent $3$-Lie algebra $(A^c,[\cdot,\cdot,\cdot]_C)$ on the vector space $V$. Then $(\rho,0)$ is a representation of the $3$-pre-Lie algebra $(A,\{\cdot,\cdot,\cdot\})$ on the vector space $V$. It is obvious that $(L,R)$ is a representation of a $3$-pre-Lie algebra on itself, which is called the {\bf regular representation}.

Let $(V,\rho,\mu)$ be a representation of a $3$-pre-Lie algebra $(A,\{\cdot,\cdot,\cdot\})$. Define a trilinear bracket operation $\{\cdot,\cdot,\cdot\}_{\rho,\mu}:\otimes^3(A\oplus V)\lon A\oplus V$  by
\begin{eqnarray}\label{semidirect-3-pre-Lie-bracket}
\{x_1+v_1,x_2+v_2,x_3+v_3\}_{\rho,\mu}\triangleq\{x_1,x_2,x_3\}+\rho(x_1,x_2)v_3+\mu(x_2,x_3)v_1-\mu(x_1,x_3)v_2.
\end{eqnarray}
By straightforward computations, we have
\begin{thm}\label{semidirect-3-pre-Lie}
With the above notation, $(A\oplus V,\{\cdot,\cdot,\cdot\}_{\rho,\mu})$ is a $3$-pre-Lie algebra.
\end{thm}
This 3-pre-Lie algebra is called the
{\bf semidirect product} of the $3$-pre-Lie algebra $(A,\{\cdot,\cdot,\cdot\})$ and $(V,\rho,\mu)$, and denoted by $A\ltimes_{\rho,\mu}V$.

Let $V$ be a vector space.   Define the switching operator $\tau:\otimes^2 V\longrightarrow \otimes^2 V$ by
\begin{eqnarray*}
\tau(T)=x_2\otimes x_1,\quad \forall T=x_1\otimes x_2\in\otimes^2 V.
\end{eqnarray*}

\begin{pro}\label{pro:representa}
Let $(\rho,\mu)$ be a representation of a $3$-pre-Lie algebra $(A,\{\cdot,\cdot,\cdot\})$ on a vector space $V$. Then $\rho-\mu\tau+\mu$ is a representation of the sub-adjacent $3$-Lie algebra $(A^c,[\cdot,\cdot,\cdot]_C)$ on the vector space $V$.
\end{pro}
\pf  By Theorem \ref{semidirect-3-pre-Lie}, we have the semidirect product 3-pre-Lie algebra $A\ltimes_{\rho,\mu}V$. Consider its sub-adjacent 3-Lie algebra structure $[\cdot,\cdot,\cdot]_C$, we have
\begin{eqnarray}
\nonumber[x_1+v_1,x_2+v_2,x_3+v_3]{_C}&=&\{x_1+v_1,x_2+v_2,x_3+v_3\}_{\rho,\mu}+\{x_2+v_2,x_3+v_3,x_1+v_1\}_{\rho,\mu}\\
                          \nonumber &&\{x_3+v_3,x_1+v_1,x_2+v_2\}_{\rho,\mu}\\
                          \nonumber &=&\{x_1,x_2,x_3\}+\rho(x_1,x_2)v_3+\mu(x_2,x_3)v_1-\mu(x_1,x_3)v_2\\
                           \nonumber&&+\{x_2,x_3,x_1\}+\rho(x_2,x_3)v_1+\mu(x_3,x_1)v_2-\mu(x_2,x_1)v_3\\
                           \nonumber&&+\{x_3,x_1,x_2\}+\rho(x_3,x_1)v_2+\mu(x_1,x_2)v_3-\mu(x_3,x_2)v_1\\
                        \nonumber   &=&[x_1,x_2,x_3]_C+((\rho-\mu\tau+\mu)(x_1,x_2))v_3\\
                          \label{eq:samesubadj} && +((\rho-\mu\tau+\mu)(x_2,x_3))v_1+((\rho-\mu\tau+\mu)(x_3,x_1))v_2.
\end{eqnarray}
By Lemma \ref{lem:semidirectp}, $\rho-\mu\tau+\mu$ is a representation of the sub-adjacent $3$-Lie algebra $(A^c,[\cdot,\cdot,\cdot]_C)$ on the vector space $V$. The proof is finished. \qed\vspace{3mm}

If $(\rho,\mu)=(L,R)$ is the regular representation of a 3-pre-Lie algebra $(A,\{\cdot,\cdot,\cdot\})$, then $\rho-\mu\tau+\mu=\ad$ is the adjoint representation of the sub-adjacent 3-Lie algebra $(A^c,[\cdot,\cdot,\cdot]_C)$.

\begin{cor}\label{sub-adjacent-3-Lie}
Let $(\rho,\mu)$ be a representation of a $3$-pre-Lie algebra $(A,\{\cdot,\cdot,\cdot\})$ on a vector space $V$. Then the semidirect product $3$-pre-Lie algebras $A\ltimes_{\rho,\mu}V$ and $A\ltimes_{\rho-\mu\tau+\mu,0}V$ given by the representations $(\rho,\mu)$ and $(\rho-\mu\tau+\mu,0)$ respectively have the same sub-adjacent $3$-Lie algebra $A^c\ltimes_{\rho-\mu\tau+\mu}V$ given by \eqref{eq:samesubadj}, which
  is the semidirect product of the $3$-Lie algebra $(A^c,[\cdot,\cdot,\cdot]_C)$ and its representation   $(V,\rho-\mu\tau+\mu)$.
\end{cor}

\emptycomment{
Let $(\rho,\mu)$ be a representation of a $3$-pre-Lie algebra $(A,\{\cdot,\cdot,\cdot\})$ on a vector space $V$. Define $\rho^*:\wedge^2 A\lon\gl(V^*)$ and $\mu^*:\otimes^2 A\lon\gl(V^*)$   by
\begin{eqnarray}
\langle \rho^*(x_1,x_2)\alpha,v\rangle=-\langle \alpha,\rho(x_1,x_2)v\rangle,\,\quad\langle \mu^*(x_1,x_2)\alpha,v\rangle=-\langle \alpha,\mu(x_1,x_2)v\rangle,
\end{eqnarray}
for all $ x_1,x_2\in A,\alpha\in V^*,v\in V.$
}

\begin{pro}
Let $(\rho,\mu)$ be a representation of a $3$-pre-Lie algebra $(A,\{\cdot,\cdot,\cdot\})$ on a vector space $V$. Then $(\rho^*-\mu^*\tau+\mu^*,-\mu^*)$ is a representation of the $3$-pre-Lie algebra $(A,\{\cdot,\cdot,\cdot\})$ on the vector space $V^*$, which is called the {\bf dual representation} of the representation $(V,\rho,\mu)$.
\end{pro}
\pf By Proposition \ref{pro:representa}, $\rho-\mu\tau+\mu$ is a representation of the sub-adjacent $3$-Lie algebra $(A^c,[\cdot,\cdot,\cdot]_C)$ on the vector space $V$. By Lemma \ref{dual-rep-3-Lie}, $\rho^*-\mu^*\tau+\mu^*$ is a representation of the sub-adjacent $3$-Lie algebra $(A^c,[\cdot,\cdot,\cdot]_C)$ on the dual vector space $V^*$. It is straightforward to deduce that other conditions in Definition \ref{defi:rep3-pre-Lie} also hold. We leave details to readers. \qed

\begin{cor}
  Let $(V,\rho,\mu)$ be a representation of a $3$-pre-Lie algebra $(A,\{\cdot,\cdot,\cdot\})$. Then the semidirect product $3$-pre-Lie algebras $A\ltimes_{\rho^*,0}V^*$ and $A\ltimes_{\rho^*-\mu^*\tau+\mu^*,-\mu^*}V^*$ given by the representations $(\rho^*,0)$ and $(\rho^*-\mu^*\tau+\mu^*,-\mu^*)$ respectively have the same sub-adjacent $3$-Lie algebra $A^c\ltimes_{\rho^*}V^*$, which
  is the semidirect product of the $3$-Lie algebra $(A^c,[\cdot,\cdot,\cdot]_C)$ and its representation $(V^*,\rho^*)$.
\end{cor}

 If $(\rho,\mu)=(L,R)$ is the regular representation of a 3-pre-Lie algebra $(A,\{\cdot,\cdot,\cdot\})$, then $(\rho^*-\mu^*\tau+\mu^*,-\mu^*)=(\ad^*,-R^*)$ and the corresponding semidirect product 3-Lie algebra is $A^c\ltimes_{L^*}A^*$, which is the key object when we construct phase spaces of 3-Lie algebras in the next section.

\section{Symplectic structures and phase spaces of $3$-Lie algebras}

In this section, we introduce the notion of a phase space of a 3-Lie algebra and show that a 3-Lie algebra has a phase space if and only if it is sub-adjacent to a 3-pre-Lie algebra. Moreover, we introduce the notion of a Manin triple of 3-pre-Lie algebras and show that there is a one-to-one correspondence between Manin triples of 3-pre-Lie algebras and perfect phase spaces of 3-Lie algebras.

\begin{defi}{\rm (\cite{BGS})}
A {\bf symplectic structure} on a $3$-Lie algebra $(\g,[\cdot,\cdot,\cdot]_\g)$   is a nondegenerate skew-symmetric bilinear form $\omega\in\wedge^2\g^*$ satisfying the following equality:
\begin{eqnarray}\label{symplectic-structure}
\omega([x,y,z]_\g,w)-\omega([y,z,w]_\g,x)+\omega([z,w,x]_\g,y)-\omega([w,x,y]_\g,z)=0,\quad\forall x,y,z,w\in\g.
\end{eqnarray}
\end{defi}

\begin{ex}\label{ex:A4symplectic}{\rm
Consider the  $4$-dimensional Euclidean $3$-Lie algebra $A_{4}$ given in \cite{BL0}. The underlying vector space is $\mathbb R^4$. Relative to an orthogonal basis $\{e_1,e_2,e_3,e_4\}$, the $3$-Lie bracket is given by
   $$[e_1,e_2,e_3]=e_4, \quad [e_2,e_3,e_4]=e_1,\quad
[e_1,e_3,e_4]=e_2,\quad[e_1,e_2,e_4]=e_3.$$
Then it is straightforward to see that any nondegenerate skew-symmetric bilinear form is a symplectic structure on $A_4$. In particular, \begin{eqnarray*}
\omega_1=e_3^*\wedge e^*_1+e_4^*\wedge e_2^*,\quad \omega_2=e_2^*\wedge e^*_1+e_4^*\wedge e_3^*,\quad \omega_3=e_2^*\wedge e^*_1+e_3^*\wedge
e_4^*,\\ \omega_4=e_1^*\wedge e^*_2+e_4^*\wedge e_3^*,\quad\omega_5=e_1^*\wedge e^*_2+e_3^*\wedge e_4^*,\quad\omega_6=e_1^*\wedge e_3^*+e_2^*\wedge e_4^*\end{eqnarray*} are symplectic structures on $A_4$, where $\{e_1^*,e_2^*,e_3^*,e_4^*\}$ are the dual basis.
}
\end{ex}

\begin{pro}{\rm (\cite{BGS})}\label{3-pre-Lie-under-3-Lie}
 Let $(\g,[\cdot,\cdot,\cdot]_\g,\omega)$ be a symplectic $3$-Lie algebra.  Then there exists
a compatible $3$-pre-Lie algebra structure $\{\cdot,\cdot,\cdot\}$ on $\g$ given by
\begin{equation}\label{3-pre-Lie-omega}
\omega(\{x,y,z\},w)=-\omega(z,[x,y,w]_\g),\quad \forall  x,y,z,w\in \g.
\end{equation}
\end{pro}

A {\bf quadratic 3-pre-Lie algebra} is a 3-pre-Lie algebra $(A,\{\cdot,\cdot,\cdot\})$ equipped with a nondegenerate skew-symmetric bilinear form $\omega\in\wedge^2A^*$ such that the following invariant condition holds:
\begin{equation}\label{eq:quadratic}
\omega(\{x,y,z\},w)=-\omega(z,[x,y,w]_C),\quad \forall  x,y,z,w\in A.
\end{equation}
Proposition \ref{3-pre-Lie-under-3-Lie} tells us that
quadratic 3-pre-Lie algebras are the underlying structures of symplectic 3-Lie algebras.

Let $V$ be a vector space and $V^*=\Hom(V,\mathbb R)$   its dual space. Then there is a natural nondegenerate skew-symmetric bilinear form $\omega$ on $T^*V=V\oplus V^*$ given by:
\begin{eqnarray}\label{phase-space}
\omega(x+\alpha,y+\beta)=\langle \alpha,y\rangle-\langle \beta,x\rangle,\,\,\,\,\forall x,y\in V,\alpha,\beta\in V^*.
\end{eqnarray}

\begin{defi}
Let $(\h,[\cdot,\cdot,\cdot]_\h)$ be a $3$-Lie algebra and $\h^*$   its dual space.
\begin{itemize}
  \item If there is a  $3$-Lie algebra structure $[\cdot,\cdot,\cdot]$ on the direct sum vector space $T^*\h=\h\oplus\h^*$ such that $(\h\oplus\h^*,[\cdot,\cdot,\cdot],\omega)$ is a symplectic $3$-Lie algebra, where $\omega$ given by \eqref{phase-space}, and $(\h,[\cdot,\cdot,\cdot]_\h)$ and $(\h^*,[\cdot,\cdot,\cdot]|_{\h^*})$ are $3$-Lie subalgebras of $ (\h\oplus\h^*,[\cdot,\cdot,\cdot])$, then the symplectic $3$-Lie algebra $(\h\oplus\h^*,[\cdot,\cdot,\cdot],\omega)$ is called a {\bf phase space} of the $3$-Lie algebra $(\h,[\cdot,\cdot,\cdot]_\h)$.

      \item A phase space $(\h\oplus\h^*,[\cdot,\cdot,\cdot],\omega)$ is called {\bf perfect} if the following conditions are satisfied:
     \begin{equation}\label{eq:conperfectPS}
      [x,y,\alpha]\in\h^*,\quad [\alpha,\beta,x]\in\h,\quad \forall x,y\in\h, \alpha,\beta\in\h^*.
     \end{equation}
\end{itemize}
\end{defi}

 3-pre-Lie algebras play important role in the study of phase spaces of 3-Lie algebras.

 \begin{thm}\label{3-pre-Lie-phase-space}
A $3$-Lie algebra has a phase space if and only if it is sub-adjacent to a $3$-pre-Lie algebra.
\end{thm}
\pf
Let $(A,\{\cdot,\cdot,\cdot\})$ be a $3$-pre-Lie algebra. By Proposition \ref{3-pre-Lie}, the left multiplication $L$ is a representation of the sub-adjacent $3$-Lie algebra $A^c$ on $A$.
By Lemma \ref{dual-rep-3-Lie}, $L^*$ is a representation of the sub-adjacent $3$-Lie algebra $A^c$ on $A^*$. Thus, we have the semidirect product 3-Lie algebra $A^c\ltimes_{L^*}A^*=(A^c\oplus A^*,[\cdot,\cdot,\cdot]_{L^*})$.
 Then $(A^c\ltimes_{L^*}A^*,\omega)$ is a symplectic $3$-Lie algebra, which is a phase space of the sub-adjacent $3$-Lie algebra $(A^c,[\cdot,\cdot,\cdot]_C)$.
In fact, for all $x_1,x_2,x_3,x_4\in A$ and $\alpha_1,\alpha_2,\alpha_3,\alpha_4\in A^*$, we have
\begin{eqnarray*}
&&\omega([x_1+\alpha_1,x_2+\alpha_2,x_3+\alpha_3]_{L^*},x_4+\alpha_4)\\&=&\omega([x_1,x_2,x_3]_C+L^*(x_1,x_2)\alpha_3+L^*(x_2,x_3)\alpha_1+L^*(x_3,x_1)\alpha_2,x_4+\alpha_4)\\
                                                       &=&\langle L^*(x_1,x_2)\alpha_3+L^*(x_2,x_3)\alpha_1+L^*(x_3,x_1)\alpha_2,x_4\rangle-\langle \alpha_4,[x_1,x_2,x_3]_C\rangle\\
                                                       &=&-\langle \alpha_3,\{x_1,x_2,x_4\}\rangle-\langle \alpha_1,\{x_2,x_3,x_4\}\rangle-\langle \alpha_2,\{x_3,x_1,x_4\}\rangle\\
                                                       &&-\langle \alpha_4,\{x_1,x_2,x_3\}\rangle-\langle \alpha_4,\{x_2,x_3,x_1\}\rangle-\langle \alpha_4,\{x_3,x_1,x_2\}\rangle.
\end{eqnarray*}
Similarly, we have
\begin{eqnarray*}
&&\omega([x_2+\alpha_2,x_3+\alpha_3,x_4+\alpha_4]_{L^*},x_1+\alpha_1)\\&=&-\langle \alpha_4,\{x_2,x_3,x_1\}\rangle-\langle \alpha_2,\{x_3,x_4,x_1\}\rangle
-\langle \alpha_3,\{x_4,x_2,x_1\}\rangle\\
                                                       &&-\langle \alpha_1,\{x_2,x_3,x_4\}\rangle-\langle \alpha_1,\{x_3,x_4,x_2\}\rangle-\langle \alpha_1,\{x_4,x_2,x_3\}\rangle,\\
&&\omega([x_3+\alpha_3,x_4+\alpha_4,x_1+\alpha_1]_{L^*},x_2+\alpha_2)\\&=&-\langle \alpha_1,\{x_3,x_4,x_2\}\rangle-\langle \alpha_3,\{x_4,x_1,x_2\}\rangle
-\langle \alpha_4,\{x_1,x_3,x_2\}\rangle\\
                                                       &&-\langle \alpha_2,\{x_3,x_4,x_1\}\rangle-\langle \alpha_2,\{x_4,x_1,x_3\}\rangle-\langle \alpha_2,\{x_1,x_3,x_4\}\rangle,\\
 &&\omega([x_4+\alpha_4,x_1+\alpha_1,x_2+\alpha_2]_{L^*},x_3+\alpha_3)\\&=&-\langle \alpha_2,\{x_4,x_1,x_3\}\rangle-\langle \alpha_4,\{x_1,x_2,x_3\}\rangle
-\langle \alpha_1,\{x_2,x_4,x_3\}\rangle\\
                                                       &&-\langle \alpha_3,\{x_4,x_1,x_2\}\rangle-\langle \alpha_3,\{x_1,x_2,x_4\}\rangle-\langle \alpha_3,\{x_2,x_4,x_1\}\rangle.
\end{eqnarray*}
Since $\{x_1,x_2,x_3\}=-\{x_2,x_1,x_3\}$, we deduce that $\omega$ is a symplectic structure on the semidirect product 3-Lie algebra $A^c\ltimes_{L^*}A^*$. Moreover, $(A^c,[\cdot,\cdot,\cdot]_C)$ is a subalgebra of $A^c\ltimes_{L^*}A^*$ and $A^*$ is an abelian subalgebra of $A^c\ltimes_{L^*}A^*$. Thus, the symplectic 3-Lie algebra $(A^c\ltimes_{L^*}A^*,\omega)$ is a phase space of the sub-adjacent $3$-Lie algebra $(A^c,[\cdot,\cdot,\cdot]_C)$.

Conversely, let $(T^*\h=\h\oplus\h^*,[\cdot,\cdot,\cdot],\omega)$ be a phase space of a $3$-Lie algebra $(\h,[\cdot,\cdot,\cdot]_\h)$. By Proposition \ref{3-pre-Lie-under-3-Lie}, there exists a compatible $3$-pre-Lie algebra structure $\{\cdot,\cdot,\cdot\}$ on $T^*\h$ given by \eqref{3-pre-Lie-omega}.
Since $(\h,[\cdot,\cdot,\cdot]_\h)$ is a subalgebra of $(\h\oplus\h^*,[\cdot,\cdot,\cdot])$,  we have
\begin{eqnarray*}
\omega(\{x,y,z\},w)=-\omega(z,[x,y,w])=-\omega(z,[x,y,w]_{\h})=0,\quad\forall x,y,z,w\in\h.
\end{eqnarray*}
Thus, $\{x,y,z\}\in\h$, which implies that $(\h,\{\cdot,\cdot,\cdot\}|_\h)$ is a subalgebra of the $3$-pre-Lie algebra $(T^*\h,\{\cdot,\cdot,\cdot\})$. Its sub-adjacent 3-Lie algebra $(\h^c,[\cdot,\cdot,\cdot]_C)$ is exactly the original $3$-Lie algebra $(\h,[\cdot,\cdot,\cdot]_\h)$.
 \qed

\begin{cor}\label{3-pre-Lie-sub}
Let $(T^*\h=\h\oplus\h^*,[\cdot,\cdot,\cdot],\omega)$ be a phase space of a $3$-Lie algebra $(\h,[\cdot,\cdot,\cdot]_\h)$ and $(\h\oplus \h^*,\{\cdot,\cdot,\cdot\})$ the associated $3$-pre-Lie algebra. Then both $(\h,\{\cdot,\cdot,\cdot\}|_\h)$ and $(\h^*,\{\cdot,\cdot,\cdot\}|_{\h^*})$ are subalgebras of the $3$-pre-Lie algebra $(\h\oplus \h^*,\{\cdot,\cdot,\cdot\})$.
\end{cor}

\begin{cor}
If $(\h\oplus\h^*,[\cdot,\cdot,\cdot],\omega)$ is a phase space of a $3$-Lie algebra $(\h,[\cdot,\cdot,\cdot]_\h)$ such that the $3$-Lie algebra $(\h\oplus\h^*,[\cdot,\cdot,\cdot])$ is a semidirect product $\h\ltimes_{\rho^*}\h^*$, where $\rho$ is a representation of $(\h,[\cdot,\cdot,\cdot]_\h)$ on  $\h$ and $\rho^*$ is its dual representation, then
$$\{x,y,z\}\triangleq \rho(x,y)z,\quad \forall x,y,z\in\h,$$
 defines a $3$-pre-Lie algebra structure on $\h$.
\end{cor}
\pf For all $x,y,z\in\h$ and $\alpha\in\h^*$, we have
\begin{eqnarray*}
\langle \alpha,\{x,y,z\}\rangle&=&-\omega(\{x,y,z\},\alpha)=\omega(z,[x,y,\alpha]_{\g\oplus\g^*})=\omega(z,\rho^*(x,y)\alpha)=-\langle \rho^*(x,y)\alpha,z\rangle\\
&=&\langle \alpha,\rho(x,y)z\rangle.
\end{eqnarray*}
Therefore, $\{x,y,z\}=\rho(x,y)z$.  \qed

\begin{ex}{\rm
Let $(A,\{\cdot,\cdot,\cdot\}_A)$ be a $3$-pre-Lie algebra. Since there is a semidirect product $3$-pre-Lie algebra structure $(A\ltimes_{L^*,0}A^*,\{\cdot,\cdot,\cdot\}_{L^*,0})$ on the phase space $T^*A^c=A^c\ltimes_{L^*}A^{*}$, one can construct a new phase space $T^*A^c\ltimes_{L^*}(T^*A^c)^*$. This process can be continued indefinitely. Hence, there exist a series of phase spaces $\{A_{(n)}\}_{n\ge2}:$
$$A_{(1)}=A^c,\,\,\,\,A_{(2)}=T^*A_{(1)}=A^c\ltimes_{L^*}A^{*},\cdots,\,\,\,\,A_{(n)}=T^*A_{(n-1)},\cdots.$$
$A_{(n)}~~(n\ge2)$ is called the symplectic double of $A_{(n-1)}.$
}
\end{ex}

At the end of this section, we introduce the notion of a Manin triple of 3-pre-Lie algebras.

\begin{defi}
  A {\bf Manin triple of $3$-pre-Lie algebras} is a triple $(\huaA;A,A')$, where
  \begin{itemize}
    \item $(\huaA,\{\cdot,\cdot,\cdot\},\omega)$ is a quadratic $3$-pre-Lie algebra;
    \item both $A$ and $A'$ are isotropic subalgebras of $(\huaA,\{\cdot,\cdot,\cdot\})$;
    \item $\huaA=A\oplus A'$ as vector spaces;
    \item for all $x,y\in A$ and $\alpha,\beta\in A'$, there holds:
    \begin{equation}\label{eq:conMT}
    \{x,y,\alpha\}\in A',\quad  \{\alpha, x,y\}\in A',\quad  \{\alpha,\beta, x\}\in A,\quad \{ x,\alpha,\beta\}\in A.
  \end{equation}
  \end{itemize}
\end{defi}

In a Manin triple of $3$-pre-Lie algebras, since the skewsymmetric bilinear form $\omega$ is nondegenerate, $A'$ can be identified with $A^*$ via
$$
\langle \alpha,x\rangle\triangleq \omega(\alpha,x),\quad\forall x\in A, \alpha\in A'.
$$
Thus, $\huaA$ is isomorphic to $A\oplus A^*$ naturally and the bilinear form $\omega$ is exactly given by \eqref{phase-space}. By the invariant condition \eqref{eq:quadratic}, we can obtain the precise form of the 3-pre-Lie structure $\{\cdot,\cdot,\cdot\}$  on $A\oplus A^*$.

\begin{pro}\label{pro:stuctureMP3preLie}
  Let $(A\oplus A^*;A,A^*)$ be a Manin triple of $3$-pre-Lie algebras, where the nondegenerate skewsymmetric bilinear form $\omega$ on the $3$-pre-Lie algebra  is  given by \eqref{phase-space}. Then we have
  \begin{eqnarray}
   \label{eq:m1}\{x,y,\alpha\}&=&(L^*-R^*\tau+R^*)(x,y)\alpha,\\
    \label{eq:m2}\{\alpha,x,y\}&=&-R^*(x,y)\alpha,\\
    \label{eq:m3}\{\alpha,\beta,x\}&=&(\huaL^*-\huaR^*\tau+\huaR^*)(\alpha,\beta)x,\\
 \label{eq:m4}\{x,\alpha,\beta\}&=&-\huaR^*(\alpha,\beta)x.
  \end{eqnarray}
\end{pro}
\pf For all $x,y,z\in A,\alpha\in A^*$, we have
\begin{eqnarray*}
\langle\{x,y,\alpha\},z\rangle&=&\omega(\{x,y,\alpha\},z)=-\omega(\alpha,[x,y,z]_C)\\
                              &=&-\omega(\alpha,\{x,y,z\}+\{y,z,x\}+\{z,x,y\})\\
                              &=&-\omega(\alpha,L(x,y)z-R(y,x)z+R(x,y)z)\\
                              &=&-\langle\alpha,L(x,y)z-R(y,x)z+R(x,y)z\rangle\\
                              &=&\langle(L^*-R^*\tau+R^*)(x,y)\alpha,z\rangle,
\end{eqnarray*}
which implies that \eqref{eq:m1} holds.  We have
\begin{eqnarray*}
\langle\{\alpha,x,y\},z\rangle&=&\omega(\{\alpha,x,y\},z)=-\omega(y,[\alpha,x,z]_C)=\omega(y,[z,x,\alpha]_C)=- \omega(\{z,x,y\},\alpha)\\
                                &=&\langle\alpha,R(x,y)z\rangle=-\langle R^*(x,y)\alpha,z\rangle,
\end{eqnarray*}
which implies that \eqref{eq:m2} holds. Similarly, we can deduce that \eqref{eq:m3} and  \eqref{eq:m4} hold.
\qed

\begin{thm}\label{thm:MT-ps}
  There is a one-to-one correspondence between Manin triples of $3$-pre-Lie algebras and perfect phase spaces of $3$-Lie algebras. More precisely, if $(A\oplus A^*;A,A^*)$ is  a Manin triple of $3$-pre-Lie algebras, then $(A\oplus A^*,[\cdot,\cdot,\cdot]_C,\omega)$ is a symplectic $3$-Lie algebra, where $\omega$ is given by \eqref{phase-space}. Conversely, if $(\h\oplus \h^*,[\cdot,\cdot,\cdot],\omega)$ is a perfect phase space of a $3$-Lie algebra $(\h,[\cdot,\cdot,\cdot]_\h)$, then $(\h\oplus \h^*;\h,\h^*)$ is  a Manin triple of $3$-pre-Lie algebras, where the $3$-pre-Lie algebra structure on $\h\oplus \h^*$ is given by \eqref{3-pre-Lie-omega}.
\end{thm}

\pf Let $(A\oplus A^*;A,A^*)$ be a Manin triple of $3$-pre-Lie algebras. Denote by $\{\cdot,\cdot,\cdot\}_A$ and $\{\cdot,\cdot,\cdot\}_{A^*}$ the 3-pre-Lie algebra structure on $A$ and $A^*$ respectively, and denote by $[\cdot,\cdot,\cdot]_A$ and $[\cdot,\cdot,\cdot]_{A^*}$ the corresponding sub-adjacent 3-Lie algebra structure on $A$ and $A^*$ respectively. By Proposition \ref{pro:stuctureMP3preLie}, it is straightforward to deduce that the corresponding 3-Lie algebra structure $[\cdot,\cdot,\cdot]_C$ on $A\oplus A^*$ is given by
\begin{eqnarray}
 \nonumber [x+\alpha,y+\beta,z+\gamma]_C&=&[x,y,z]_A+\huaL^*(\alpha,\beta)z+\huaL^*(\beta,\gamma)x+\huaL^*(\gamma,\alpha)y\\
\label{eq:MP3Lie} &&+[\alpha,\beta,\gamma]_{A^*}+L^*(x,y)\gamma+L^*(y,z)\alpha+L^*(z,x)\beta.
\end{eqnarray}
 For all $x_1,x_2,x_3,x_4\in A$ and $\alpha_1,\alpha_2,\alpha_3,\alpha_4\in A^*$, we have
\begin{eqnarray*}
&&\omega([x_1+\alpha_1,x_2+\alpha_2,x_3+\alpha_3]_C,x_4+\alpha_4)\\&=&\omega([x_1,x_2,x_3]_A+\huaL^*(\alpha_1,\alpha_2)x_3+\huaL^*(\alpha_2,\alpha_3)x_1+\huaL^*(\alpha_3,\alpha_1)x_2\\
                     &&+[\alpha_1,\alpha_2,\alpha_3]_{A^*}+L^*(x_1,x_2)\alpha_3+L^*(x_2,x_3)\alpha_1+L^*(x_3,x_1)\alpha_2,x_4+\alpha_4)\\
                     &=&\langle [\alpha_1,\alpha_2,\alpha_3]_{A^*}+L^*(x_1,x_2)\alpha_3+L^*(x_2,x_3)\alpha_1+L^*(x_3,x_1)\alpha_2,x_4\rangle\\
                     &&-\langle \alpha_4,[x_1,x_2,x_3]_A+\huaL^*(\alpha_1,\alpha_2)x_3+\huaL^*(\alpha_2,\alpha_3)x_1+\huaL^*(\alpha_3,\alpha_1)x_2\rangle\\
                     &=&\langle [\alpha_1,\alpha_2,\alpha_3]_{A^*},x_4\rangle-\langle \alpha_3,\{x_1,x_2,x_4\}_A\rangle-\langle \alpha_1,\{x_2,x_3,x_4\}_A\rangle-\langle \alpha_2,\{x_3,x_1,x_4\}_A\rangle\\
                     &&-\langle \alpha_4,[x_1,x_2,x_3]_A\rangle+\langle\{\alpha_1,\alpha_2,\alpha_4\}_{A^*},x_3\rangle+\langle\{\alpha_2,\alpha_3,\alpha_4\}_{A^*},x_1\rangle
                     +\langle\{\alpha_3,\alpha_1,\alpha_4\}_{A^*},x_2\rangle.
\end{eqnarray*}
Similarly, we have
 \begin{eqnarray*}
&&\omega([x_2+\alpha_2,x_3+\alpha_3,x_4+\alpha_4],x_1+\alpha_1)\\&=&\langle [\alpha_2,\alpha_3,\alpha_4]_C,x_1\rangle-\langle \alpha_4,\{x_2,x_3,x_1\}_A\rangle-\langle \alpha_2,\{x_3,x_4,x_1\}_A\rangle-\langle \alpha_3,\{x_4,x_2,x_1\}_A\rangle\\
                     &&-\langle \alpha_1,[x_2,x_3,x_4]_C\rangle+\langle\{\alpha_2,\alpha_3,\alpha_1\}_{A^*},x_4\rangle+\langle\{\alpha_3,\alpha_4,\alpha_1\}_{A^*},x_2\rangle
                     +\langle\{\alpha_4,\alpha_2,\alpha_1\}_{A^*},x_3\rangle,\\
&&\omega([x_3+\alpha_3,x_4+\alpha_4,x_1+\alpha_1],x_2+\alpha_2)\\&=&\langle [\alpha_3,\alpha_4,\alpha_1]_C,x_2\rangle-\langle \alpha_1,\{x_3,x_4,x_2\}_A\rangle-\langle \alpha_3,\{x_4,x_1,x_2\}_A\rangle-\langle \alpha_4,\{x_1,x_3,x_2\}_A\rangle\\
                     &&-\langle \alpha_2,[x_3,x_4,x_1]_C\rangle+\langle\{\alpha_3,\alpha_4,\alpha_2\}_{A^*},x_1\rangle+\langle\{\alpha_4,\alpha_1,\alpha_2\}_{A^*},x_3\rangle
                     +\langle\{\alpha_1,\alpha_3,\alpha_2\}_{A^*},x_4\rangle,\\
&&\omega([x_4+\alpha_4,x_1+\alpha_1,x_2+\alpha_2],x_3+\alpha_3)\\&=&\langle [\alpha_4,\alpha_1,\alpha_2]_C,x_3\rangle-\langle \alpha_2,\{x_4,x_1,x_3\}_A\rangle-\langle \alpha_4,\{x_1,x_2,x_3\}_A\rangle-\langle \alpha_1,\{x_2,x_4,x_3\}_A\rangle\\
                     &&-\langle \alpha_3,[x_4,x_1,x_2]_C\rangle+\langle\{\alpha_4,\alpha_1,\alpha_3\}_{A^*},x_2\rangle+\langle\{\alpha_1,\alpha_2,\alpha_3\}_{A^*},x_4\rangle
                     +\langle\{\alpha_2,\alpha_4,\alpha_3\}_{A^*},x_1\rangle.
\end{eqnarray*}
By $\{x_1,x_2,x_3\}_A=-\{x_2,x_1,x_3\}_A$ and $\{\alpha_1,\alpha_2,\alpha_3\}_{A^*}=-\{\alpha_2,\alpha_1,\alpha_3\}_{A^*}$, we deduce that $\omega$ is a symplectic structure on the 3-Lie algebra $(A\oplus A^*,[\cdot,\cdot,\cdot]_C)$. Therefore, it is a phase space.

Conversely, let $(\h\oplus \h^*,[\cdot,\cdot,\cdot],\omega)$ be a phase space of the $3$-Lie algebra $(\h,[\cdot,\cdot,\cdot]_\h)$. By Proposition \ref{3-pre-Lie-under-3-Lie}, there exists a $3$-pre-Lie algebra structure $ \{\cdot,\cdot,\cdot\}$ on $ \h\oplus \h^*$ given by  \eqref{3-pre-Lie-omega} such that $(\h\oplus \h^*,\{\cdot,\cdot,\cdot\},\omega)$ is a quadratic 3-pre-Lie algebra. By Corollary \ref{3-pre-Lie-sub}, $(\h,\{\cdot,\cdot,\cdot\}|_{\h})$ and $(\h^*,\{\cdot,\cdot,\cdot\}|_{\h^*})$ are $3$-pre-Lie subalgebras of $(\h\oplus \h^*,\{\cdot,\cdot,\cdot\})$. It is obvious that both $\h$ and $\h^*$ are isotropic.  Thus, we only need to show that \eqref{eq:conMT} holds. By \eqref{eq:conperfectPS}, for all $x_1,x_2\in\h$ and $\alpha_1,\alpha_2\in\h^*$,  we have
\begin{eqnarray*}
  \omega(\{x_1,x_2,\alpha_1\},\alpha_2)=-\omega(\alpha_1,[x_1,x_2,\alpha_2]_C)=0,
\end{eqnarray*}
which implies that $\{x_1,x_2,\alpha_1\}\in\h^*$.  Similarly, we can show that the other conditions in \eqref{eq:conMT} also hold.  The proof is finished. \qed

\begin{rmk}
  The notions of a matched pair of $3$-Lie algebras and a Manin triple of $3$-Lie algebras were introduced in \cite{BGS}. By \eqref{eq:MP3Lie}, we obtain that $(A^c,{A^*}^c;L^*,\huaL^*)$ is a matched pair of $3$-Lie algebras and the phase space is exactly the double of this matched pair. However, one should note that a Manin triple of $3$-pre-Lie algebras does not give rise to a Manin triple of $3$-Lie algebras.
\end{rmk}
\begin{rmk}
For pre-Lie algebras, there are equivalent description between Manin triples of pre-Lie algebras, matched pairs of pre-Lie algebras associated to the dual representations of the regular representations and pre-Lie bialgebras \cite{Baibialgebra}.   Here we only study Manin triples of $3$-pre-Lie algebras, which are closely related to phase spaces of $3$-Lie algebras and para-Kähler $3$-Lie algebras that studied in Section 8.  We postpone the study of matched pairs of $3$-pre-Lie algebras and $3$-pre-Lie bialgebras in the future.
\end{rmk}

\section{Product structures on $3$-Lie algebras}

In this section, we introduce the notion of a product structure on a 3-Lie algebra using the Nijenhuis condition as the integrability condition. We find four special integrability conditions, each of them gives a special decomposition of the original 3-Lie algebra. At the end of this section, we introduce the notion of a (perfect) paracomplex structure on a $3$-Lie algebra and give examples.

\begin{defi}
Let $(\g,[\cdot,\cdot,\cdot]_\g)$ be a $3$-Lie algebra. An {\bf almost product structure} on the $3$-Lie algebra $(\g,[\cdot,\cdot,\cdot]_\g)$ is a linear endomorphism $E:\g\lon\g$ satisfying $E^2=\Id$ $(E\not=\pm\Id)$. An  almost product structure   is called a {\bf product} structure if the  following integrability   condition is satisfied:
\begin{eqnarray}\nonumber\label{product-structure}
E[x,y,z]_\g&=&[Ex,Ey,Ez]_\g+[Ex,y,z]_\g+[x,Ey,z]_\g+[x,y,Ez]_\g\\
&&-E[Ex,Ey,z]_\g-E[x,Ey,Ez]_\g-E[Ex,y,Ez]_\g.
\end{eqnarray}
\end{defi}

 \begin{rmk}
   One can understand a product structure on a $3$-Lie algebra as a Nijenhuis operator $E$ on a $3$-Lie algebra satisfying $E^2=\Id.$
 \end{rmk}

\begin{thm}\label{product-structure-subalgebra}
Let $(\g,[\cdot,\cdot,\cdot]_\g)$ be a $3$-Lie algebra. Then $(\g,[\cdot,\cdot,\cdot]_\g)$ has a product structure if and only if $\g$ admits a decomposition:
\begin{eqnarray}
\g=\g_+\oplus\g_-,
\end{eqnarray}
where $\g_+$ and $\g_-$ are subalgebras of $\g$.
\end{thm}
\pf Let $E$ be a product structure on $\g$. By $E^2=\Id$ , we have $\g=\g_+\oplus\g_-$, where $\g_+$ and $\g_-$ are the eigenspaces of $\g$ associated to the eigenvalues $\pm1$. For all $x_1,x_2,x_3\in\g_+$, we have
\begin{eqnarray*}
E[x_1,x_2,x_3]_\g&=&[Ex_1,Ex_2,Ex_3]_\g+[Ex_1,x_2,x_3]_\g+[x_1,Ex_2,x_3]_\g+[x_1,x_2,Ex_3]_\g\\
&&-E[Ex_1,Ex_2,x_3]_\g-E[x_1,Ex_2,Ex_3]_\g-E[Ex_1,x_2,Ex_3]_\g\\
&=&4[x_1,x_2,x_3]_\g-3E[x_1,x_2,x_3]_\g.
\end{eqnarray*}
Thus, we have $[x_1,x_2,x_3]_\g\in\g_{+}$, which implies that $\g_+$ is a subalgebra. Similarly, we can show that $\g_-$ is a subalgebra.

Conversely, we define a linear endomorphism $E:\g\lon\g$ by
\begin{eqnarray}\label{eq:productE}
E(x+\alpha)=x-\alpha,\,\,\,\,\forall x\in\g_+,\alpha\in\g_-.
\end{eqnarray}
Obviously we have $E^2=\Id$.
Since $\g_+$ is a subalgebra of $\g$, for all $x_1,x_2,x_3\in\g_+$, we have
\begin{eqnarray*}
&&[Ex_1,Ex_2,Ex_3]_\g+[Ex_1,x_2,x_3]_\g+[x_1,Ex_2,x_3]_\g+[x_1,x_2,Ex_3]_\g\\
&&-E[Ex_1,Ex_2,x_3]_\g-E[x_1,Ex_2,Ex_3]_\g-E[Ex_1,x_2,Ex_3]_\g\\
&=&4[x_1,x_2,x_3]_\g-3E[x_1,x_2,x_3]_\g=[x_1,x_2,x_3]_\g\\
&=&E[x_1,x_2,x_3]_\g,
\end{eqnarray*}
which implies that \eqref{product-structure} holds for all $x_1,x_2,x_3\in\g_+$. Similarly, we can show that \eqref{product-structure} holds for all $x,y,z\in\g$.
Therefore,   $E$ is a product structure on $\g$.  \qed

\begin{lem}
Let $E$ be an almost product structure on a $3$-Lie algebra $(\g,[\cdot,\cdot,\cdot]_\g)$. If $E$ satisfies the following equation
\begin{eqnarray}\label{abel-product-0}
E[x,y,z]_\g=[Ex,y,z]_\g,
\end{eqnarray}
then $E$ is a product structure on $\g$ such that $[\g_+,\g_+,\g_-]_\g=0$ and $[\g_-,\g_-,\g_+]_\g=0,$ i.e. $\g$ is the $3$-Lie algebra direct sum of $\g_+$ and $\g_-.$
\end{lem}
\pf By \eqref{abel-product-0} and $E^2=\Id$, we have
\begin{eqnarray*}
&&[Ex,Ey,Ez]_\g+[Ex,y,z]_\g+[x,Ey,z]_\g+[x,y,Ez]_\g\\
&&-E[Ex,Ey,z]_\g-E[x,Ey,Ez]_\g-E[Ex,y,Ez]_\g\\
&=&[Ex,Ey,Ez]_\g+E[x,y,z]_\g+[x,Ey,z]_\g+[x,y,Ez]_\g\\
&&-[E^2x,Ey,z]_\g-[Ex,Ey,Ez]_\g-[E^2x,y,Ez]_\g\\
&=&E[x,y,z]_\g.
\end{eqnarray*}
Thus, $E$ is a product structure on $\g$. For all $x_1,x_2\in\g_+,\alpha_1\in\g_-$, on one hand we have
\begin{eqnarray*}
E[\alpha_1,x_1,x_2]_\g=[E\alpha_1,x_1,x_2]_\g=-[\alpha_1,x_1,x_2]_\g.
\end{eqnarray*}
On the other hand, we have
\begin{eqnarray*}
E[\alpha_1,x_1,x_2]_\g=E[x_1,x_2,\alpha_1]_\g=[Ex_1,x_2,\alpha_1]_\g=[x_1,x_2,\alpha_1]_\g.
\end{eqnarray*}
Thus, we obtain $[\g_+,\g_+,\g_-]_\g=0$. Similarly,   we have
  $[\g_-,\g_-,\g_+]_\g=0$. The proof is finished. \qed

\begin{defi}{\bf (Integrability condition I)}
An almost product structure  $E$ on  a $3$-Lie algebra  $(\g,[\cdot,\cdot,\cdot]_\g)$ is called  a {\bf strict product structure} if  \eqref{abel-product-0} holds.
\end{defi}

\begin{cor}
Let $(\g,[\cdot,\cdot,\cdot]_\g)$ be a $3$-Lie algebra. Then $(\g,[\cdot,\cdot,\cdot]_\g)$ has a strict product structure if and only if
$\g$ admits a decomposition:
$$
\g=\g_+\oplus\g_-,
$$
where $\g_+$ and $\g_-$ are subalgebras of $\g$ such that $[\g_+,\g_+,\g_-]_\g=0$ and $[\g_-,\g_-,\g_+]_\g=0.$
\end{cor}
\pf We leave the details to readers. \qed

\begin{lem}
Let $E$ be an almost product structure on a $3$-Lie algebra $(\g,[\cdot,\cdot,\cdot]_\g)$. If $E$ satisfies the following equation
\begin{eqnarray}\label{abel-product}
[x,y,z]_\g=-[x,Ey,Ez]_\g-[Ex,y,Ez]_\g-[Ex,Ey,z]_\g,
\end{eqnarray}
then $E$ is a product structure on $\g$.
\end{lem}
\pf By \eqref{abel-product} and $E^2=\Id$, we have
\begin{eqnarray*}
&&[Ex,Ey,Ez]_\g+[Ex,y,z]_\g+[x,Ey,z]_\g+[x,y,Ez]_\g\\
&&-E[Ex,Ey,z]_\g-E[x,Ey,Ez]_\g-E[Ex,y,Ez]_\g\\
&=&-[Ex,E^2y,E^2z]_\g-[E^2x,Ey,E^2z]_\g-[E^2x,E^2y,Ez]_\g\\
&&+[Ex,y,z]_\g+[x,Ey,z]_\g+[x,y,Ez]_\g+E[x,y,z]_\g\\
&=&E[x,y,z]_\g.
\end{eqnarray*}
Thus,   $E$ is a product structure on $\g$. \qed

\begin{defi}{\bf (Integrability condition II)}
An almost product structure  $E$ on a $3$-Lie algebra  $(\g,[\cdot,\cdot,\cdot]_\g)$ is called  an {\bf abelian product structure} if  \eqref{abel-product} holds.
\end{defi}

\begin{cor}\label{abelian-product-structure}
Let $(\g,[\cdot,\cdot,\cdot]_\g)$ be a $3$-Lie algebra. Then $(\g,[\cdot,\cdot,\cdot]_\g)$ has an abelian  product structure if and only if
$\g$ admits a decomposition:
$$
\g=\g_+\oplus\g_-,
$$
where $\g_+$ and $\g_-$ are abelian subalgebras of $\g$.
\end{cor}
\pf Let $E$ be an abelian product structure on $\g$. For all $x_1,x_2,x_3\in\g_{+}$, we have
\begin{eqnarray*}
[x_1,x_2,x_3]_\g&=& -[Ex_1,Ex_2,x_3]_\g-[x_1,Ex_2,Ex_3]_\g-[Ex_1,x_2,Ex_3]_\g\\
           &=&-3[x_1,x_2,x_3]_\g,
\end{eqnarray*}
which implies that $[x_1,x_2,x_3]_\g=0$. Similarly, for all $\alpha_1,\alpha_2,\alpha_3\in\g_{-}$, we also have
  $[\alpha_1,\alpha_2,\alpha_3]_\g=0$. Thus, both $\g_+$ and $\g_-$ are abelian subalgebras.

Conversely,   define a linear endomorphism $E:\g\lon\g$ by \eqref{eq:productE}.
Then it is straightforward to deduce that $E$ is an abelian product structure on $\g$.   \qed

\begin{lem}
Let $E$ be an almost product structure on a $3$-Lie algebra $(\g,[\cdot,\cdot,\cdot]_\g)$. If $E$ satisfies the following equation
\begin{eqnarray}\label{product-integrability}
[x,y,z]_\g&=&E[Ex,y,z]_\g+E[x,Ey,z]_\g+E[x,y,Ez]_\g,
\end{eqnarray}
then $E$ is an abelian product structure on $\g$ such that $[\g_+,\g_+,\g_-]_\g\subset\g_+$ and $[\g_-,\g_-,\g_+]_\g\subset\g_-.$
\end{lem}
\pf By \eqref{product-integrability} and $E^2=\Id$, we have
\begin{eqnarray*}
&&[Ex,Ey,Ez]_\g+[Ex,y,z]_\g+[x,Ey,z]_\g+[x,y,Ez]_\g\\
&&-E[Ex,Ey,z]_\g-E[x,Ey,Ez]_\g-E[Ex,y,Ez]_\g\\
&=&E[x,Ey,Ez]_\g+E[Ex,y,Ez]_\g+E[Ex,Ey,z]_\g+E[x,y,z]_\g\\
&&-E[Ex,Ey,z]_\g-E[x,Ey,Ez]_\g-E[Ex,y,Ez]_\g\\
&=&E[x,y,z]_\g.
\end{eqnarray*}
Thus, we obtain that $E$ is a product structure on $\g$. For all $x_1,x_2,x_3\in\g_+$, by \eqref{product-integrability}, we have
\begin{eqnarray*}
[x_1,x_2,x_3]_\g&=&E[Ex_1,x_2,x_3]_\g+E[x_1,Ex_2,x_3]_\g+E[x_1,x_2,Ex_3]_\g\\
                &=&3E[x_1,x_2,x_3]_\g=3[x_1,x_2,x_3]_\g.
\end{eqnarray*}
Thus, we obtain $[\g_+,\g_+,\g_+]_\g=0$. Similarly,   we have $[\g_-,\g_-,\g_-]_\g=0$. By Corollary \ref{abelian-product-structure}, $E$ is an abelian product structure on $\g$. Moreover, for all $x_1,x_2\in\g_+,\alpha_1\in\g_-$, we have
\begin{eqnarray*}
[x_1,x_2,\alpha_1]_\g&=&E[Ex_1,x_2,\alpha_1]_\g+E[x_1,Ex_2,\alpha_1]_\g+E[x_1,x_2,E\alpha_1]_\g\\
                &=&E[x_1,x_2,\alpha_1]_\g,
\end{eqnarray*}
which implies that $[\g_+,\g_+,\g_-]_\g\subset\g_+$. Similarly,    we have
  $[\g_-,\g_-,\g_+]_\g\subset\g_-$.   \qed

\begin{defi}{\bf (Integrability condition III)}
An almost product structure $E$ on a $3$-Lie algebra $ (\g,[\cdot,\cdot,\cdot]_\g)$  is called a {\bf strong abelian product structure} if   \eqref{product-integrability} holds.
\end{defi}

\begin{cor}
Let $(\g,[\cdot,\cdot,\cdot]_\g)$ be a $3$-Lie algebra. Then $(\g,[\cdot,\cdot,\cdot]_\g)$ has a strong abelian product structure if and only if
$\g$ admits a decomposition:
$$
\g=\g_+\oplus\g_-,
$$
where $\g_+$ and $\g_-$ are abelian subalgebras of $\g$ such that $[\g_+,\g_+,\g_-]_\g\subset\g_+$ and $[\g_-,\g_-,\g_+]_\g\subset\g_-.$
\end{cor}
\begin{rmk}
 Let $E$ be a strong abelian product structure on a $3$-Lie algebra $(\g,[\cdot,\cdot,\cdot]_\g)$. Then we can define $\nu_+:\g_+\longrightarrow \Hom(\wedge^2\g_-,\g_-)$ and $\nu_-:\g_-\longrightarrow \Hom(\wedge^2\g_+,\g_+)$ by
 $$
 \nu_+(x)(\alpha,\beta)=[\alpha,\beta,x]_\g,\quad \nu_-(\alpha)(x,y)=[x,y,\alpha]_\g,\quad\forall x,y\in\g_+, \alpha,\beta\in\g_-.
 $$
 It turns out $\nu_+$ and $\nu_-$ are generalized representations of abelian $3$-Lie algebras $\g_+$ and $\g_-$ on $\g_-$ and $\g_+$ respectively. See \cite{Liu-Makhlouf-Sheng} for more details about generalized representations of $3$-Lie algebras.
\end{rmk}

More surprisingly, a strong abelian product structure is an $\huaO$-operator as well as a Rota-Baxter operator \cite{RB3Lie,PBG}. Thus, some $\huaO$-operators and Rota-Baxter operators on 3-Lie algebras can serve as integrability conditions.
\begin{pro}
Let $E$ be an almost product structure on a $3$-Lie algebra $ (\g,[\cdot,\cdot,\cdot]_\g)$. Then $E$ is a strong abelian structure on $\g$ if and only if $E$ is an $\huaO$-operator associated to the adjoint representation $(\g,\ad)$. Furthermore, there exists a compatible $3$-pre-Lie algebra $(\g,\{\cdot,\cdot,\cdot\})$ on the $3$-Lie algebra $(\g,[\cdot,\cdot,\cdot]_\g)$, here the $3$-pre-Lie algebra structure on $\g$ is given by
\begin{eqnarray}
\{x,y,z\}=E[x,y,Ez]_\g,\,\,\,\,\forall x,y,z\in\g.
\end{eqnarray}
\end{pro}
\pf By \eqref{product-integrability}, for all $x,y,z\in\g$ we have
\begin{eqnarray*}
[Ex,Ey,Ez]_\g&=&E[E^2x,Ey,Ez]_\g+E[Ex,E^2y,Ez]_\g+E[Ex,Ey,E^2z]_\g\\
             &=&E(\ad_{Ex,Ey}z+\ad_{Ey,Ez}x+\ad_{Ez,Ex}y).
\end{eqnarray*}
Thus, $E$ is an $\huaO$-operator associated to the adjoint representation $(\g,\ad)$.

Conversely, if for  all $x,y,z\in\g$, we have
\begin{eqnarray*}
[Ex,Ey,Ez]_\g&=&E(\ad_{Ex,Ey}z+\ad_{Ey,Ez}x+\ad_{Ez,Ex}y)\\
             &=&E([Ex,Ey,z]_\g+[x,Ey,Ez]_\g+[Ex,y,Ez]_\g),
\end{eqnarray*}
 then
$
[x,y,z]_\g=E[x,y,Ez]_\g+E[Ex,y,z]_\g+E[x,Ey,z]_\g
$ by   $E^{-1}=E$.
Thus, $E$ is a strong abelian structure on $\g$.

Furthermore, by $E^{-1}=E$ and Proposition \ref{3-Lie-compatible-3-pre-Lie}, there exists a compatible $3$-pre-Lie algebra on $\g$ given by
$
\{x,y,z\}=E\ad_{x,y}E^{-1}(z)=E[x,y,Ez]_\g.
$
The proof is finished. \qed\vspace{3mm}

There is a new phenomenon that an involutive automorphism of a 3-Lie algebra also serves as an integrability condition.

\begin{lem}
Let $E$ be an almost product structure on a $3$-Lie algebra $(\g,[\cdot,\cdot,\cdot]_\g)$. If $E$ satisfies the following equation
\begin{eqnarray}\label{product-integrability-1}
E[x,y,z]_\g=[Ex,Ey,Ez]_\g,
\end{eqnarray}
then $E$ is a product structure on $\g$ such that \begin{equation}\label{eq:coherenceconPP}
[\g_+,\g_+,\g_-]_\g\subset\g_-,\quad [\g_-,\g_-,\g_+]_\g\subset\g_+.
\end{equation}
\end{lem}
\pf By \eqref{product-integrability-1} and $E^2=\Id$, we have
\begin{eqnarray*}
&&[Ex,Ey,Ez]_\g+[Ex,y,z]_\g+[x,Ey,z]_\g+[x,y,Ez]_\g\\
&&-E[Ex,Ey,z]_\g-E[x,Ey,Ez]_\g-E[Ex,y,Ez]_\g\\
&=&E[x,y,z]_\g+[Ex,y,z]_\g+[x,Ey,z]_\g+[x,y,Ez]_\g\\
&&-[E^2x,E^2y,Ez]_\g-[Ex,E^2y,E^2z]_\g-[E^2x,Ey,E^2z]_\g\\
&=&E[x,y,z]_\g.
\end{eqnarray*}
Thus,   $E$ is a product structure on $\g$. Moreover, for all $x_1,x_2\in\g_+,\alpha_1\in\g_-$, we have
\begin{eqnarray*}
E[x_1,x_2,\alpha_1]_\g=[Ex_1,Ex_2,E\alpha_1]_\g=-[x_1,x_2,\alpha_1]_\g,
\end{eqnarray*}
which implies that $[\g_+,\g_+,\g_-]_\g\subset\g_-$. Similarly, we have $[\g_-,\g_-,\g_+]_\g\subset\g_+$.  \qed

\begin{defi}{\bf (Integrability condition IV)}
An almost product structure $E$ on a $3$-Lie algebra $(\g,[\cdot,\cdot,\cdot]_\g)$ is  called   a {\bf perfect product structure} if   \eqref{product-integrability-1} holds.
\end{defi}

\begin{cor}
Let $(\g,[\cdot,\cdot,\cdot]_\g)$ be a $3$-Lie algebra. Then $(\g,[\cdot,\cdot,\cdot]_\g)$ has a perfect product structure if and only if
$\g$ admits a decomposition:
$$
\g=\g_+\oplus\g_-,
$$
where $\g_+$ and $\g_-$ are subalgebras of $\g$ such that $[\g_+,\g_+,\g_-]_\g\subset\g_-$ and $[\g_-,\g_-,\g_+]_\g\subset\g_+.$
\end{cor}
\pf We leave the details to readers. \qed

\begin{cor}
  A strict product structure on a $3$-Lie algebra    is a perfect product structure.
\end{cor}
\begin{rmk}
Let $E$ be a   product structure on a $3$-Lie algebra $(\g,[\cdot,\cdot,\cdot]_\g)$. By Theorem \ref{product-structure-subalgebra}, $\g_+$ and $\g_-$ are subalgebras. However, the brackets of mixed terms are very complicated. But a perfect product structure $E$   on $(\g,[\cdot,\cdot,\cdot]_\g)$ ensures   $[\g_+,\g_+,\g_-]_\g\subset\g_-$ and $[\g_-,\g_-,\g_+]_\g\subset\g_+$. Note that this is exactly the condition required in the definition of a matched pair of $3$-Lie algebras \cite{BGS}. Thus, $E$ is a perfect product structure if and only if $(\g_+,\g_-)$ is a matched pair of $3$-Lie algebras. This type of product structures are very important in our later studies.
\end{rmk}

\emptycomment{
\begin{eqnarray*}
[E(x_1+y_1),E(x_2+y_2),E(x_3+y_3)]_\g&=&[x_1-y_1,x_2-y_2,x_3-y_3]_\g\\
                                     &=&[x_1,x_2,-y_3]_\g+[x_1,-y_2,x_3]_\g+[x_1,-y_2,-y_3]_\g\\
                                     &&+[-y_1,x_2,x_3]_\g+[-y_1,x_2,-y_3]_\g+[-y_1,-y_2,x_3]_\g,\\
\,[E(x_1+y_1),x_2+y_2,x_3+y_3]_\g&=&[x_1-y_1,x_2+y_2,x_3+y_3]_\g\\
                                 &=&[x_1,x_2,y_3]_\g+[x_1,y_2,x_3]_\g+[x_1,y_2,y_3]_\g\\
                                 &&+[-y_1,x_2,x_3]_\g+[-y_1,x_2,y_3]_\g+[-y_1,y_2,x_3]_\g,\\
\,[x_1+y_1,E(x_2+y_2),x_3+y_3]_\g&=&[x_1+y_1,x_2-y_2,x_3+y_3]_\g\\
                                 &=&[x_1,x_2,y_3]_\g+[x_1,-y_2,x_3]_\g+[x_1,-y_2,y_3]_\g\\
                                     &&+[y_1,x_2,x_3]_\g+[y_1,x_2,y_3]_\g+[y_1,-y_2,x_3]_\g,\\
\,[x_1+y_1,x_2+y_2,E(x_3+y_3)]_\g&=&[x_1+y_1,x_2+y_2,x_3-y_3]_\g\\
                                 &=&[x_1,x_2,-y_3]_\g+[x_1,y_2,x_3]_\g+[x_1,y_2,-y_3]_\g\\
                                     &&+[y_1,x_2,x_3]_\g+[y_1,x_2,-y_3]_\g+[y_1,y_2,x_3]_\g.
\end{eqnarray*}
By \eqref{product-structure}, we obtain $$E[x,y,z]_\g=-E[Ex_1,Ex_2,x_3]_\g-E[x_1,Ex_2,Ex_3]_\g-E[Ex_1,x_2,Ex_3]_\g.$$
Thus, we obtain $E$ is an abelian product structure on $\g$. The proof is finished. \qed
}

\begin{defi}
\begin{itemize}
  \item[\rm(i)]A {\bf paracomplex structure} on a $3$-Lie algebra $(\g,[\cdot,\cdot,\cdot]_\g)$ is a product structure $E$ on $\g$ such that the eigenspaces of $\g$ associated to the eigenvalues $\pm1$  have the same dimension, i.e. $\dim(\g_+)=\dim(\g_-)$.

    \item[\rm(i)]A {\bf perfect paracomplex structure} on a $3$-Lie algebra $(\g,[\cdot,\cdot,\cdot]_\g)$ is a perfect product structure $E$ on $\g$ such that the eigenspaces of $\g$ associated to the eigenvalues $\pm1$  have the same dimension, i.e. $\dim(\g_+)=\dim(\g_-)$.
\end{itemize}

\end{defi}

\begin{pro}\label{paracomplex-3-pre-Lie}
Let $(A,\{\cdot,\cdot,\cdot\})$ be a $3$-pre-Lie algebra. Then, on the semidirect product $3$-Lie algebra $ A^c\ltimes_{L^*}A^*$, there is a perfect paracomplex structure $E:A^c\ltimes_{L^*}A^*\lon A^c\ltimes_{L^*}A^*$ given by
\begin{eqnarray}\label{eq:defiE}
E(x+\alpha)=x-\alpha,\,\,\,\,\forall x\in A^c, \alpha\in A^*.
\end{eqnarray}
\end{pro}
\pf It is obvious that $E^2=\Id$. Moreover, we have $(A^c\ltimes_{L^*}A^*)_+=A$, $(A^c\ltimes_{L^*}A^*)_-=A^*$
and they are two subalgebras of the semidirect product $3$-Lie algebra $A^c\ltimes_{L^*}A^*$. By Theorem \ref{product-structure-subalgebra},   $E$ is a product structure on $A^c\ltimes_{L^*}A^*$. Since $A$ and $A^*$ have the same dimension, $E$ is a paracomplex structure on $A^c\ltimes_{L^*}A^*$.  It is obvious that $E$ is perfect. \qed\vspace{3mm}

At the end of this section,  we give some examples of product structures.
\begin{ex}{\rm
 There is a unique non-trivial $3$-dimensional $3$-Lie algebra. It has a basis $\{e_1,e_2,e_3\}$ with respect to which the non-zero product is given by
$$[e_1,e_2,e_3]=e_1.$$
Then $E=\left(\begin{array}{ccc}1&0&0\\
0&1&0\\
0&0&-1\end{array}\right)$ and $E=\left(\begin{array}{ccc}1&0&0\\
0&-1&0\\
0&0&1\end{array}\right)$ are strong abelian product structures and $E=\left(\begin{array}{ccc}-1&0&0\\
0&1&0\\
0&0&1\end{array}\right)$ is a perfect product structure.
}
\end{ex}

\begin{ex}\label{ex:A4product}{\rm
Consider the $4$-dimensional Euclidean $3$-Lie algebra $A_4$ given in Example \ref{ex:A4symplectic}.
Then \begin{eqnarray*}E_1=\left(\begin{array}{cccc}1&0&0&0\\
0&1&0&0\\
0&0&-1&0\\
0&0&0&-1\end{array}\right),~
E_2=\left(\begin{array}{cccc}1&0&0&0\\
0&-1&0&0\\
0&0&1&0\\
0&0&0&-1\end{array}\right),~
E_3=\left(\begin{array}{cccc}1&0&0&0\\
0&-1&0&0\\
0&0&-1&0\\
0&0&0&1\end{array}\right),\\
E_4=\left(\begin{array}{cccc}-1&0&0&0\\
0&1&0&0\\
0&0&1&0\\
0&0&0&-1\end{array}\right),~
E_5=\left(\begin{array}{cccc}-1&0&0&0\\
0&1&0&0\\
0&0&-1&0\\
0&0&0&1\end{array}\right),
E_6=\left(\begin{array}{cccc}-1&0&0&0\\
0&-1&0&0\\
0&0&1&0\\
0&0&0&1\end{array}\right)
\end{eqnarray*}are perfect and abelian product structures.
}
\end{ex}

\section{Complex structures on $3$-Lie algebras}

In this section, we introduce the notion of a complex structure on a real 3-Lie algebra using the Nijenhuis condition as the integrability condition. Parallel to the case of product structures, we also find four special integrability conditions.
\begin{defi}\label{complex}
Let $(\g,[\cdot,\cdot,\cdot]_\g)$ be a real $3$-Lie algebra. An {\bf almost complex structure} on $\g$ is a linear endomorphism $J:\g\lon\g$ satisfying $J^2=-\Id$. An   almost complex structure is called a  {\bf complex} structure if the following integrability  condition is satisfied:
\begin{eqnarray}\nonumber\label{complex-structure}
J[x,y,z]_\g&=&-[Jx,Jy,Jz]_\g+[Jx,y,z]_\g+[x,Jy,z]_\g+[x,y,Jz]_\g\\
&&+J[Jx,Jy,z]_\g+J[x,Jy,Jz]_\g+J[Jx,y,Jz]_\g.
\end{eqnarray}
 \end{defi}

 \begin{rmk}
   One can understand a complex structure on a $3$-Lie algebra as a Nijenhuis operator $J$ on a $3$-Lie algebra satisfying $J^2=-\Id.$
 \end{rmk}

\begin{rmk}
One can also use definition \ref{complex} to define the notion of a complex structure on a
complex $3$-Lie algebra, considering $J$ to be $\mathbb C$-linear. However, this is not very interesting since
for a complex $3$-Lie algebra, there is a one-to-one correspondence between such $\mathbb C$-linear complex
structures and product structures  (see Proposition \ref{equivalent}).
\end{rmk}

Consider  $\g_{\mathbb C}=\g\otimes_{\mathbb R} \mathbb C\cong\{x+iy|x,y\in\g\}$, the complexification of the
real Lie algebra $\g$, which turns out to be a complex $3$-Lie algebra by extending the $3$-Lie bracket on $\g$ complex trilinearly, and we  denote it by $(\g_{\mathbb C},[\cdot,\cdot,\cdot]_{\g_{\mathbb C}})$. We have an equivalent description of the integrability condition given in
Definition \ref{complex}. We denote by $\sigma$ the conjugation in $\g_{\mathbb C}$ with respect to the real form $\g$,
that is, $\sigma(x+iy)=x-iy,\,\,x,y\in\g$. Then, $\sigma$ is a complex antilinear, involutive automorphism of the complex vector space $\g_{\mathbb C}$.

\begin{thm}\label{complex-structure-subalgebra}
Let $(\g,[\cdot,\cdot,\cdot]_\g)$ be a real $3$-Lie algebra. Then $\g$ has a complex structure if and only if $\g_{\mathbb C}$ admits a decomposition:
\begin{eqnarray}
\g_{\mathbb C}=\frkq\oplus\frkp,
\end{eqnarray}
where $\frkq$ and $\frkp=\sigma(\frkq)$ are complex subalgebras of $\g_{\mathbb C}$.
\end{thm}
\pf We  extend the complex structure $J$ complex linearly, which is denoted by $J_{\mathbb C}$, i.e. $J_{\mathbb C}:\g_{\mathbb C}\longrightarrow \g_{\mathbb C}$ is defined by
\begin{equation}\label{eq:JC}
J_{\mathbb C}(x+iy)=Jx+iJy,\quad \forall x,y\in\g.
\end{equation} Then $J_{\mathbb C}$ is a complex linear endomorphism on $\g_{\mathbb C}$ satisfying $J_{\mathbb C}^2=-\Id$ and the integrability condition \eqref{complex-structure} on $\g_{\mathbb C}$. Denote by $\g_{\pm i}$ the corresponding eigenspaces of $\g_{\mathbb C}$ associated to the eigenvalues $\pm i$ and there holds:
   \begin{eqnarray*}
\g_{\mathbb C}=\g_{i}\oplus\g_{-i}.
\end{eqnarray*}
  It is straightforward to see that  $\g_{i}=\{x-iJx|x\in\g\}$ and $\g_{-i}=\{x+iJx|x\in\g\}$. Therefore, we have $\g_{-i}=\sigma(\g_{i})$.

For all $X,Y,Z\in\g_{i}$, we have
\begin{eqnarray*}
J_{\mathbb C}[X,Y,Z]_{\g_{\mathbb C}}&=&-[J_{\mathbb C}X,J_{\mathbb C}Y,J_{\mathbb C}Z]_{\g_{\mathbb C}}+[J_{\mathbb C}X,Y,Z]_{\g_{\mathbb C}}+[X,J_{\mathbb C}Y,Z]_{\g_{\mathbb C}}+[X,Y,J_{\mathbb C}Z]_{\g_{\mathbb C}}\\
&&+J_{\mathbb C}[J_{\mathbb C}X,J_{\mathbb C}Y,Z]_{\g_{\mathbb C}}+J_{\mathbb C}[X,J_{\mathbb C}Y,J_{\mathbb C}Z]_{\g_{\mathbb C}}+J_{\mathbb C}[J_{\mathbb C}X,Y,J_{\mathbb C}Z]_{\g_{\mathbb C}}\\
&=&4i[X,Y,Z]_{\g_{\mathbb C}}-3J_{\mathbb C}[X,Y,Z]_{\g_{\mathbb C}}.
\end{eqnarray*}
Thus, we have $[X,Y,Z]_{\g_{\mathbb C}}\in\g_{i}$, which implies that $\g_i$ is a subalgebra. Similarly, we can show that $\g_{-i}$ is also a subalgebra.

Conversely, we define a complex linear endomorphism $J_{\mathbb C}:\g_{\mathbb C}\lon\g_{\mathbb C}$ by
\begin{eqnarray}\label{defi-complex-structure}
J_{\mathbb C}(X+\sigma(Y))=iX-i\sigma(Y),\,\,\,\,\forall X,Y\in\frkq.
\end{eqnarray}
Since $\sigma$ is a complex antilinear, involutive automorphism of $\g_{\mathbb C}$, we have
\begin{eqnarray*}
J_{\mathbb C}^2(X+\sigma(Y))=J_{\mathbb C}(iX-i\sigma(Y))=J_{\mathbb C}(iX+\sigma(iY))=i(iX)-i\sigma(iY)=-X-\sigma(Y),
\end{eqnarray*}
i.e. $J_{\mathbb C}^2=-\Id$. Since $\frkq$ is a subalgebra of $\g_{\mathbb C}$, for all $X,Y,Z\in\frkq$, we have
\begin{eqnarray*}
&&-[J_{\mathbb C}X,J_{\mathbb C}Y,J_{\mathbb C}Z]_{\g_{\mathbb C}}+[J_{\mathbb C}X,Y,Z]_{\g_{\mathbb C}}+[X,J_{\mathbb C}Y,Z]_{\g_{\mathbb C}}+[X,Y,J_{\mathbb C}Z]_{\g_{\mathbb C}}\\
&&+J_{\mathbb C}[J_{\mathbb C}X,J_{\mathbb C}Y,Z]_{\g_{\mathbb C}}+J_{\mathbb C}[X,J_{\mathbb C}Y,J_{\mathbb C}Z]_{\g_{\mathbb C}}+J_{\mathbb C}[J_{\mathbb C}X,Y,J_{\mathbb C}Z]_{\g_{\mathbb C}}\\
&=&4i[X,Y,Z]_{\g_{\mathbb C}}-3J_{\mathbb C}[X,Y,Z]_{\g_{\mathbb C}}=i[X,Y,Z]_{\g_{\mathbb C}}\\
&=&J_{\mathbb C}[X,Y,Z]_{\g_{\mathbb C}},
\end{eqnarray*}
which implies that $J_{\mathbb C}$ satisfies   \eqref{complex-structure} for all $X,Y,Z\in\frkq$. Similarly, we can show that $J_{\mathbb C}$ satisfies   \eqref{complex-structure} for all $\huaX,\huaY,\huaZ\in\g_{\mathbb C}$. Since $\g_{\mathbb C}=\frkq\oplus\frkp$, we can write $\huaX\in\g_{\mathbb C}$   as $\huaX=X+\sigma(Y),$ for some $X,Y\in\frkq$. Since $\sigma$ is a complex antilinear, involutive automorphism of $\g_{\mathbb C}$, we have
\begin{eqnarray*}
(J_{\mathbb C}\circ\sigma)(X+\sigma(Y))=J_{\mathbb C}(Y+\sigma(X))=iY-i\sigma(X)=\sigma(iX-i\sigma(Y))=(\sigma\circ J_{\mathbb C})(X+\sigma(Y)),
\end{eqnarray*}
which implies that $J_{\mathbb C}\circ\sigma=\sigma\circ J_{\mathbb C}$. Moreover, since $\sigma(\huaX)=\huaX$ is equivalent to $\huaX\in\g$, we deduce that the set of fixed points of $\sigma$ is the real vector
space $\g$. By $J_{\mathbb C}\circ\sigma=\sigma\circ J_{\mathbb C}$, there is a well-defined $J\in\gl(\g)$ given by $$J\triangleq J_{\mathbb C}|_{\g}.$$ Follows from that $J_{\mathbb C}$  satisfies \eqref{complex-structure} and $J_{\mathbb C}^2=-\Id$ on $\g_{\mathbb C}$, $J$ is a complex structure on $\g$.  \qed

\begin{lem}
Let $J$ be an almost complex structure on a real $3$-Lie algebra $(\g,[\cdot,\cdot,\cdot]_\g)$. If $J$ satisfies
\begin{eqnarray}\label{adapt}
J[x,y,z]_\g=[Jx,y,z]_\g,\,\,\,\,\forall x,y,z\in\g,
\end{eqnarray}
  then $J$ is a complex structure on $(\g,[\cdot,\cdot,\cdot]_\g)$.
\end{lem}

\pf By \eqref{adapt} and $J^2=-\Id$, we have
\begin{eqnarray*}
&&-[Jx,Jy,Jz]_\g+[Jx,y,z]_\g+[x,Jy,z]_\g+[x,y,Jz]_\g\\
&&+J[Jx,Jy,z]_\g+J[x,Jy,Jz]_\g+J[Jx,y,Jz]_\g\\
&=&-[Jx,Jy,Jz]_\g+J[x,y,z]_\g+[x,Jy,z]_\g+[x,y,Jz]_\g\\
&&+[J^2x,Jy,z]_\g+[Jx,Jy,Jz]_\g+[J^2x,y,Jz]_\g\\
&=&J[x,y,z]_\g.
\end{eqnarray*}
Thus, we obtain that $J$ is a complex structure on $\g$. \qed

\begin{defi}{\bf (Integrability condition I)}
  An almost complex structure $J$  on a real $3$-Lie algebra $(\g,[\cdot,\cdot,\cdot]_\g)$ is called a {\bf strict complex structure} if \eqref{adapt} holds.
\end{defi}

\begin{cor}
Let $(\g,[\cdot,\cdot,\cdot]_\g)$ be a real $3$-Lie algebra. Then there is a strict  complex structure on $(\g,[\cdot,\cdot,\cdot]_\g)$ if and only if $\g_{\mathbb C}$ admits a decomposition:
\begin{eqnarray}
\g_{\mathbb C}=\frkq\oplus\frkp,
\end{eqnarray}
where $\frkq$ and $\frkp=\sigma(\frkq)$ are complex subalgebras of $\g_{\mathbb C}$ such that $[\frkq,\frkq,\frkp]_{\g_{\mathbb C}}=0$ and $[\frkp,\frkp,\frkq]_{\g_{\mathbb C}}=0,$ i.e. $\g_{\mathbb C}$ is a $3$-Lie algebra direct sum of $\frkq$ and $\frkp$.
\end{cor}
\pf
Let $J$ be a strict complex structure on a real $3$-Lie algebra $(\g,[\cdot,\cdot,\cdot]_\g)$. Then, $J_{\mathbb C}$ is a strict complex structure on the complex $3$-Lie algebra $(\g_{\mathbb C},[\cdot,\cdot,\cdot]_{\g_{\mathbb C}})$.   For all $X,Y\in\g_{i}$ and $\sigma(Z)\in\g_{-i}$, on one hand we have
\begin{eqnarray*}
J_{\mathbb C}[X,Y,\sigma(Z)]_{\g_{\mathbb C}}=[J_{\mathbb C}X,Y,\sigma(Z)]_{\g_{\mathbb C}}=i[X,Y,\sigma(Z)]_{\g_{\mathbb C}}.
\end{eqnarray*}
On the other hand, we have
\begin{eqnarray*}
J_{\mathbb C}[X,Y,\sigma(Z)]_{\g_{\mathbb C}}=J_{\mathbb C}[\sigma(Z),X,Y]_{\g_{\mathbb C}}=[J_{\mathbb C}\sigma(Z),X,Y]_{\g_{\mathbb C}}=-i[\sigma(Z),X,Y]_{\g_{\mathbb C}}.
\end{eqnarray*}
Thus, we obtain $[\g_i,\g_i,\g_{-i}]_{\g_{\mathbb C}}=0$. Similarly, we can show $[\g_{-i},\g_{-i},\g_{i}]_{\g_{\mathbb C}}=0$.

Conversely, define a complex linear endomorphism $J_{\mathbb C}:\g_{\mathbb C}\lon\g_{\mathbb C}$ by \eqref{defi-complex-structure}. Then it is straightforward to deduce that $J_{\mathbb C}^2=-\Id$. Since $\frkq$ is a subalgebra of $\g_{\mathbb C}$, for all $X,Y,Z\in\frkq$, we have
\begin{eqnarray*}
J_{\mathbb C}[X,Y,Z]_{\g_{\mathbb C}}=i[X,Y,Z]_{\g_{\mathbb C}}=[J_{\mathbb C}X,Y,Z]_{\g_{\mathbb C}},
\end{eqnarray*}
which implies that $J_{\mathbb C}$ satisfies   \eqref{adapt} for all $X,Y,Z\in\frkq$. Similarly, we can show that $J_{\mathbb C}$ satisfies   \eqref{adapt} for all $\huaX,\huaY,\huaZ\in\g_{\mathbb C}$. By the proof of Theorem \ref{complex-structure-subalgebra}, we obtain that $J\triangleq J_{\mathbb C}|_{\g}$ is a  strict complex structure on   the real $3$-Lie algebra $(\g,[\cdot,\cdot,\cdot]_\g)$. The proof is finished. \qed\vspace{2mm}

Let $J$ be an almost complex structure on a real 3-Lie algebra $(\g,[\cdot,\cdot,\cdot]_\g)$. We can define a complex vector
space structure on the real vector space $\g$ by
\begin{eqnarray}\label{complex-space}
(a+bi)x\triangleq ax+bJx,\,\,\,\forall a,b\in\mathbb R,x\in\g.
\end{eqnarray}
Define two maps $\varphi:\g\lon\g_i$ and $\psi:\g\lon\g_{-i}$ as following:
\begin{eqnarray*}
\varphi(x)&=&\frac{1}{2}(x-iJx),\\
\psi(x)   &=&\frac{1}{2}(x+iJx).
\end{eqnarray*}
It is straightforward to deduce that $\varphi$ is complex linear isomorphism and $\psi=\sigma\circ\varphi$ is a complex antilinear isomorphism between complex vector spaces.

Let $J$ be a strict complex structure on a real 3-Lie algebra  $(\g,[\cdot,\cdot,\cdot]_\g)$. Then with the complex vector space structure  defined above,  $(\g,[\cdot,\cdot,\cdot]_\g)$ is a complex $3$-Lie algebra.
 In fact,  the fact that the $3$-Lie bracket is complex trilinear follows from
\begin{eqnarray*}
[(a+bi)x,y,z]_\g&=&[ax+bJx,y,z]_\g=a[x,y,z]_\g+b[Jx,y,z]_\g\\
                &=&a[x,y,z]_\g+bJ[x,y,z]_\g=(a+bi)[x,y,z]_\g
\end{eqnarray*}
using \eqref{adapt} and \eqref{complex-space}. \vspace{2mm}

Let $J$ be a complex structure on $\g$. Define a new bracket $[\cdot,\cdot,\cdot]_J:\wedge^3\g\lon\g$ by
\begin{eqnarray}\label{J-bracket}
[x,y,z]_J\triangleq \frac{1}{4}([x,y,z]_\g-[x,Jy,Jz]_\g-[Jx,y,Jz]_\g-[Jx,Jy,z]_\g),\,\,\,\,\forall x,y,z\in\g.
\end{eqnarray}

\begin{pro}\label{subalgebra-iso}
Let $J$ be a complex structure on a real $3$-Lie algebra $(\g,[\cdot,\cdot,\cdot]_\g)$. Then $(\g,[\cdot,\cdot,\cdot]_J)$ is a real $3$-Lie algebra. Moreover, $J$ is a strict complex structure on  $(\g,[\cdot,\cdot,\cdot]_J)$ and the corresponding complex $3$-Lie algebra $(\g,[\cdot,\cdot,\cdot]_J)$  is isomorphic to the complex $3$-Lie algebra $\g_{i}$.
\end{pro}
\pf One can show that $(\g,[\cdot,\cdot,\cdot]_J)$ is a real $3$-Lie algebra directly. Here we use a different approach to prove this result. By \eqref{complex-structure}, for all $x,y,z\in\g$, we have
\begin{eqnarray}
\nonumber[\varphi(x),\varphi(y),\varphi(z)]_{\g_{\mathbb C}}&=&\frac{1}{8}[x-iJx,y-iJy,z-iJz]_{\g_{\mathbb C}}\\
                                                  \nonumber &=&\frac{1}{8}([x,y,z]_\g-[x,Jy,Jz]_\g-[Jx,y,Jz]_\g-[Jx,Jy,z]_\g)\\
                                                   \nonumber&&-\frac{1}{8}i([x,y,Jz]_\g+[x,Jy,z]_\g+[Jx,y,z]_\g-[Jx,Jy,Jz]_\g)\\
                                                  \nonumber &=&\frac{1}{8}([x,y,z]_\g-[x,Jy,Jz]_\g-[Jx,y,Jz]_\g-[Jx,Jy,z]_\g)\\
                                                  \nonumber &&-\frac{1}{8}iJ([x,y,z]_\g-[x,Jy,Jz]_\g-[Jx,y,Jz]_\g-[Jx,Jy,z]_\g)\\
                                                   \label{eq:Jiso}&=&\varphi[x,y,z]_J.
\end{eqnarray}
Thus, we have $[x,y,z]_J=\varphi^{-1}[\varphi(x),\varphi(y),\varphi(z)]_{\g_{\mathbb C}}$. Since $J$ is a complex structure, $\g_i$ is a $3$-Lie subalgebra.   Therefore, $(\g,[\cdot,\cdot,\cdot]_J)$ is a real $3$-Lie algebra.

By \eqref{complex-structure}, for all $x,y,z\in\g$, we have
\begin{eqnarray*}
J[x,y,z]_J&=&\frac{1}{4}J([x,y,z]_\g-[x,Jy,Jz]_\g-[Jx,y,Jz]_\g-[Jx,Jy,z]_\g)\\
          &=&\frac{1}{4}(-[Jx,Jy,Jz]_\g+[Jx,y,z]_\g+[x,Jy,z]_\g+[x,y,Jz]_\g)\\
          &=&[Jx,y,z]_J,
\end{eqnarray*}
which implies that $J$ is a strict complex structure on  $(\g,[\cdot,\cdot,\cdot]_J)$. By \eqref{eq:Jiso},  $\varphi$ is a complex $3$-Lie algebra isomorphism. The proof is finished. \qed

\begin{pro}
Let $J$ be a complex structure on a real $3$-Lie algebra $(\g,[\cdot,\cdot,\cdot]_\g)$. Then $J$ is a strict complex structure on  $(\g,[\cdot,\cdot,\cdot]_\g)$ if and only if $[\cdot,\cdot,\cdot]_J=[\cdot,\cdot,\cdot]_\g.$
\end{pro}
\pf If $J$ is a strict complex structure on $(\g,[\cdot,\cdot,\cdot]_\g)$, by $J[x,y,z]_\g=[Jx,y,z]_\g$, we have
\begin{eqnarray*}
[x,y,z]_J=\frac{1}{4}([x,y,z]_\g-[x,Jy,Jz]_\g-[Jx,y,Jz]_\g-[Jx,Jy,z]_\g)=[x,y,z]_\g.
\end{eqnarray*}
Conversely, if $[\cdot,\cdot,\cdot]_J=[\cdot,\cdot,\cdot]_\g$, we have
$$-3[x,y,z]_\g=[x,Jy,Jz]_\g+[Jx,y,Jz]_\g+[Jx,Jy,z]_\g.$$
Then by the integrability condition of $J$, we obtain
\begin{eqnarray*}
4J[x,y,z]_J&=&-[Jx,Jy,Jz]_\g+[Jx,y,z]_\g+[x,Jy,z]_\g+[x,y,Jz]_\g\\
            &=&3[Jx,y,z]_\g+[Jx,y,z]_\g\\
            &=&4[Jx,y,z]_\g,
\end{eqnarray*}
which implies that $J[x,y,z]_\g=[Jx,y,z]_\g$. The proof is finished. \qed

\begin{lem}
Let $J$ be an almost complex structure on a real $3$-Lie algebra $(\g,[\cdot,\cdot,\cdot]_\g)$. If $J$ satisfies the following equation
\begin{eqnarray}\label{abel-complex}
[x,y,z]_\g=[x,Jy,Jz]_\g+[Jx,y,Jz]_\g+[Jx,Jy,z]_\g,
\end{eqnarray}
Then, $J$ is a complex structure on $\g$.
\end{lem}
\pf By \eqref{abel-complex} and $J^2=-\Id$, we have
\begin{eqnarray*}
&&-[Jx,Jy,Jz]_\g+[Jx,y,z]_\g+[x,Jy,z]_\g+[x,y,Jz]_\g\\
&&+J[Jx,Jy,z]_\g+J[x,Jy,Jz]_\g+J[Jx,y,Jz]_\g\\
&=&-[Jx,J^2y,J^2z]_\g-[J^2x,Jy,J^2z]_\g-[J^2x,J^2y,Jz]_\g\\
&&+[Jx,y,z]_\g+[x,Jy,z]_\g+[x,y,Jz]_\g+J[x,y,z]_\g\\
&=&J[x,y,z]_\g.
\end{eqnarray*}
Thus, we obtain that $J$ is a complex structure on $\g$. \qed

\begin{defi}{\bf (Integrability condition II)}
An almost complex structure $J$ on a real $3$-Lie algebra $(\g,[\cdot,\cdot,\cdot]_\g)$ is called an {\bf abelian complex structure} if  \eqref{abel-complex} holds.
\end{defi}

\begin{rmk}
Let $J$ be an abelian complex structure on a real $3$-Lie algebra $(\g,[\cdot,\cdot,\cdot]_\g)$. Then $(\g,[\cdot,\cdot,\cdot]_J)$ is an abelian $3$-Lie algebra.
\end{rmk}

\begin{cor}
Let $(\g,[\cdot,\cdot,\cdot]_\g)$ be a real $3$-Lie algebra. Then $\g$ has an abelian complex structure if and only if $\g_{\mathbb C}$ admits a decomposition:
$$
\g_{\mathbb C}=\frkq\oplus\frkp,
$$
where $\frkq$ and $\frkp=\sigma(\frkq)$ are complex abelian subalgebras of $\g_{\mathbb C}$.
\end{cor}
\pf Let $J$ be an abelian complex structure on $\g$. By Proposition \ref{subalgebra-iso}, we obtain that $\varphi$ is a complex $3$-Lie algebra isomorphism from $(\g,[\cdot,\cdot,\cdot]_J)$ to $(\g_{i},[\cdot,\cdot,\cdot]_{\g_{\mathbb C}})$. Since $J$ is abelian, $(\g,[\cdot,\cdot,\cdot]_J)$ is an abelian $3$-Lie algebra. Therefore, $\frkq=\g_{i}$ is an abelian subalgebra of $\g_{\mathbb C}$. Since $\frkp=\g_{-i}=\sigma(\g_{i})$, for all $x_1+iy_1,x_2+iy_2,x_3+iy_3\in\g_i$, we have
\begin{eqnarray*}
&&[\sigma(x_1+iy_1),\sigma(x_2+iy_2),\sigma(x_3+iy_3)]_{\g_{\mathbb C}}\\&=&[x_1-iy_1,x_2-iy_2,x_3-iy_3]_{\g_{\mathbb C}}\\
                                                                     &=&([x_1,x_2,x_3]_\g-[x_1,y_2,y_3]_\g-[y_1,x_2,y_3]_\g-[y_1,y_2,x_3]_\g)\\
                                                                     &&-i([x_1,x_2,y_3]_\g+[x_1,y_2,x_3]_\g+[y_1,x_2,x_3]_\g-[y_1,y_2,y_3]_\g)\\
                                                                     &=&\sigma[x_1+iy_1,x_2+iy_2,x_3+iy_3]_{\g_{\mathbb C}}\\
                                                                     &=&0.
\end{eqnarray*}
Thus, $\frkp$ is an abelian subalgebra of $\g_{\mathbb C}$.

Conversely, by Theorem \ref{complex-structure-subalgebra}, there is a complex structure $J$ on $\g$. Moreover, by Proposition \ref{subalgebra-iso}, we have a complex $3$-Lie algebra isomorphism $\varphi$ from $(\g,[\cdot,\cdot,\cdot]_J)$ to $(\frkq,[\cdot,\cdot,\cdot]_{\g_{\mathbb C}})$. Thus, $(\g,[\cdot,\cdot,\cdot]_J)$ is an abelian $3$-Lie algebra. By the definition of $[\cdot,\cdot,\cdot]_J$, we obtain that $J$ is an abelian complex structure on $\g$. The proof is finished. \qed

\begin{lem}
Let $J$ be an almost complex structure on a real $3$-Lie algebra $(\g,[\cdot,\cdot,\cdot]_\g)$. If $J$ satisfies the following equation
\begin{eqnarray}\label{complex-integrability}
[x,y,z]_\g&=&-J[Jx,y,z]_\g-J[x,Jy,z]_\g-J[x,y,Jz]_\g,
\end{eqnarray}
then $J$ is a complex structure on $\g$.
\end{lem}
\pf By \eqref{complex-integrability} and $J^2=-\Id$, we have
\begin{eqnarray*}
&&-[Jx,Jy,Jz]_\g+[Jx,y,z]_\g+[x,Jy,z]_\g+[x,y,Jz]_\g\\
&&+J[Jx,Jy,z]_\g+J[x,Jy,Jz]_\g+J[Jx,y,Jz]_\g\\
&=&J[J^2x,Jy,Jz]_\g+J[Jx,J^2y,Jz]_\g+J[Jx,Jy,J^2z]_\g+J[x,y,z]_\g\\
&&+J[Jx,Jy,z]_\g+J[x,Jy,Jz]_\g+J[Jx,y,Jz]_\g\\
&=&J[x,y,z]_\g.
\end{eqnarray*}
Thus, $J$ is a complex structure on $\g$.   \qed

\begin{defi}{\bf (Integrability condition III)}
An almost complex structure $J$   on a real $3$-Lie algebra $(\g,[\cdot,\cdot,\cdot]_\g)$ is called a {\bf strong abelian complex structure} if  \eqref{complex-integrability} holds.
\end{defi}

\begin{cor}
Let $(\g,[\cdot,\cdot,\cdot]_\g)$ be a real $3$-Lie algebra. Then $\g$ has a strong abelian complex structure if and only if $\g_{\mathbb C}$ admits a decomposition:
$$
\g_{\mathbb C}=\frkq\oplus\frkp,
$$
where $\frkq$ and $\frkp=\sigma(\frkq)$ are abelian complex subalgebras of $\g_{\mathbb C}$ such that $[\frkq,\frkq,\frkp]_{\g_{\mathbb C}}\subset\frkq$ and $[\frkp,\frkp,\frkq]_{\g_{\mathbb C}}\subset\frkp.$
\end{cor}

Parallel to the case of strong abelian product structures on a 3-Lie algebra, strong abelian complex structures on a $3$-Lie algebra are also $\huaO$-operators associated to the adjoint representation.

\begin{pro}
Let $J$ be an almost complex structure on a real $3$-Lie algebra $(\g,[\cdot,\cdot,\cdot]_\g)$. Then $J$ is a strong abelian complex structure on a $3$-Lie algebra $ (\g,[\cdot,\cdot,\cdot]_\g)$ if and only if  $-J$ is an $\huaO$-operator on $(\g,[\cdot,\cdot,\cdot]_\g)$ associated to the adjoint representation $(\g,\ad)$. Furthermore, there exists a compatible $3$-pre-Lie algebra $(\g,\{\cdot,\cdot,\cdot\})$ on the $3$-Lie algebra $(\g,[\cdot,\cdot,\cdot]_\g)$, here the $3$-pre-Lie algebra structure on $\g$ is given by
\begin{eqnarray}
\{x,y,z\}=-J[x,y,Jz]_\g,\,\,\,\,\forall x,y,z\in\g.
\end{eqnarray}
\end{pro}
\pf By \eqref{complex-integrability}, for all $x,y,z\in\g$ we have
\begin{eqnarray*}
[-Jx,-Jy,-Jz]_\g&=&J[J^2x,Jy,Jz]_\g+J[Jx,J^2y,Jz]_\g+J[Jx,Jy,J^2z]_\g\\
             &=&-J(\ad_{-Jx,-Jy}z+\ad_{-Jy,-Jz}x+\ad_{-Jz,-Jx}y).
\end{eqnarray*}
Thus, $-J$ is an $\huaO$-operator associated to the adjoint representation $(\g,\ad)$.

Conversely, if for all $x,y,z\in\g$, we have
\begin{eqnarray*}
[-Jx,-Jy,-Jz]_\g&=&-J(\ad_{-Jx,-Jy}z+\ad_{-Jy,-Jz}x+\ad_{-Jz,-Jx}y)\\
             &=&-J([-Jx,-Jy,z]_\g+[x,-Jy,-Jz]_\g+[-Jx,y,-Jz]_\g),
\end{eqnarray*}
 then we obtain
$[x,y,z]_\g=-J[x,y,Jz]_\g-J[Jx,y,z]_\g-J[x,Jy,z]_\g$ by $({-J})^{-1}=J$.

Furthermore, by $(-J)^{-1}=J$ and Proposition \ref{3-Lie-compatible-3-pre-Lie}, there exists a compatible $3$-pre-Lie algebra on $\g$ given by
$
\{x,y,z\}=-J\ad_{x,y}({-J}^{-1}(z))=-J[x,y,Jz]_\g.
$
The proof is finished. \qed

\begin{lem}
Let $J$ be an almost complex structure on a real $3$-Lie algebra $(\g,[\cdot,\cdot,\cdot]_\g)$. If $J$ satisfies the following equation
\begin{eqnarray}\label{complex-integrability-1}
J[x,y,z]_\g=-[Jx,Jy,Jz]_\g,
\end{eqnarray}
then $J$ is a complex structure on $\g$.
\end{lem}
\pf By \eqref{complex-integrability-1} and $J^2=\Id$, we have
\begin{eqnarray*}
&&-[Jx,Jy,Jz]_\g+[Jx,y,z]_\g+[x,Jy,z]_\g+[x,y,Jz]_\g\\
&&+J[Jx,Jy,z]_\g+J[x,Jy,Jz]_\g+J[Jx,y,Jz]_\g\\
&=&J[x,y,z]_\g+[Jx,y,z]_\g+[x,Jy,z]_\g+[x,y,Jz]_\g\\
&&-[J^2x,J^2y,Jz]_\g-[Jx,J^2y,J^2z]_\g-[J^2x,Jy,J^2z]_\g\\
&=&J[x,y,z]_\g.
\end{eqnarray*}
Thus, $J$ is a complex structure on $\g$.  \qed

\begin{defi}{\bf (Integrability condition IV)}
An almost complex structure $J$ on a real $3$-Lie algebra $(\g,[\cdot,\cdot,\cdot]_\g)$ is called a {\bf perfect complex structure} if  \eqref{complex-integrability-1} holds.
\end{defi}

\begin{cor}
Let $(\g,[\cdot,\cdot,\cdot]_\g)$ be a real $3$-Lie algebra. Then $\g$ has a perfect complex structure if and only if $\g_{\mathbb C}$ admits a decomposition:
$$
\g_{\mathbb C}=\frkq\oplus\frkp,
$$
where $\frkq$ and $\frkp=\sigma(\frkq)$ are complex subalgebras of $\g_{\mathbb C}$ such that $[\frkq,\frkq,\frkp]_{\g_{\mathbb C}}\subset\frkp$ and $[\frkp,\frkp,\frkq]_{\g_{\mathbb C}}\subset\frkq.$
\end{cor}

\begin{cor}
Let $J$ be a strict complex structure on a real $3$-Lie algebra $(\g,[\cdot,\cdot,\cdot]_\g)$. Then $J$ is a perfect complex structure on $\g$.
\end{cor}

\begin{ex}\label{ex:A4complex}{\rm
Consider the $4$-dimensional Euclidean $3$-Lie algebra $A_4$ given in Example \ref{ex:A4symplectic}. Then
\begin{eqnarray*}
J_1=\left(\begin{array}{cccc}0&0&-1&0\\
                            0&0&0&-1\\
                            1&0&0&0\\
                            0&1&0&0\end{array}\right),~
J_2=\left(\begin{array}{cccc}0&-1&0&0\\
                            1&0&0&0\\
                            0&0&0&-1\\
                            0&0&1&0\end{array}\right),~
J_3=\left(\begin{array}{cccc}0&-1&0&0\\
                           1&0&0&0\\
                           0&0&0&1\\
                           0&0&-1&0\end{array}\right),~\\
J_4=\left(\begin{array}{cccc} 0&1&0&0\\
                           -1&0&0&0\\
                            0&0&0&-1\\
                            0&0&1&0\end{array}\right),~
J_5=\left(\begin{array}{cccc}  0&1&0&0\\
                            -1&0&0&0\\
                             0&0&0&1\\
                             0&0&-1&0\end{array}\right),~
J_6=\left(\begin{array}{cccc}  0&0&1&0\\
                             0&0&0&1\\
                            -1&0&0&0\\
                             0&-1&0&0\end{array}\right)\end{eqnarray*}are abelian complex structures. Moreover, $J_1,J_6$ are strong abelian complex structures and $J_2,J_3,J_4,J_5$ are perfect complex structures.

}
\end{ex}

\section{Complex product structures on $3$-Lie algebras}

In this section, we add a compatibility condition between a complex structure and a product structure on a 3-Lie algebra to introduce the notion of a complex product structure. We construct complex product structures using 3-pre-Lie algebras. First we illustrate the relation between a complex structure and a product structure on a complex $3$-Lie algebra.

\begin{pro}\label{equivalent}
Let $(\g,[\cdot,\cdot,\cdot]_\g)$ be a complex $3$-Lie algebra. Then $E$ is a product structure on $\g$ if and only if $J=iE$ is a complex structure on $\g$.
\end{pro}

\pf Let $E$ be a product structure on $\g$. We have $J^2=i^2E^2=-\Id.$ Thus, $J$ is an almost complex structure on $\g$. Since $E$ satisfies the integrability condition \eqref{product-structure}, we have
\begin{eqnarray*}
J[x,y,z]_\g&=&iE[x,y,z]_\g\\
                 &=&-[iEx,iEy,iEz]_\g+[iEx,y,z]_\g+[x,iEy,z]_\g+[x,y,iEz]_\g\\
                 &&+iE[iEx,iEy,z]_\g+iE[x,iEy,iEz]_\g+iE[iEx,y,iEz]_\g.
\end{eqnarray*}
Thus, $J$ is a complex structure on the complex $3$-Lie algebra $\g$.

The converse part can be proved similarly and we omit details. \qed

\begin{cor}\label{complex-to-special-paracomplex}
Let $J$ be a complex structure on a real $3$-Lie algebra $(\g,[\cdot,\cdot,\cdot]_\g)$. Then, $-iJ_{\mathbb C}$ is a paracomplex structure on the complex $3$-Lie algebra $(\g_{\mathbb C},[\cdot,\cdot,\cdot]_{\g_{\mathbb C}})$, where $J_{\mathbb C}$ is defined by \eqref{eq:JC}.
\end{cor}
\pf By Theorem \ref{complex-structure-subalgebra},   $\g_{\mathbb C}=\g_{i}\oplus\g_{-i}$ and $\g_{-i}=\sigma(\g_{i})$,
where $\g_{i}$ and $\g_{-i}$ are subalgebras of $\g_{\mathbb C}$. It is obvious that $\dim(\g_i)=\dim(\g_{-i})$.  By Proposition \ref{product-structure-subalgebra}, there is a paracomplex structure on $\g_{\mathbb C}$. On the other hand, it is obvious that $J_{\mathbb C}$ is a complex structure on $\g_{\mathbb C}$. By Proposition \ref{equivalent}, $-iJ_{\mathbb C}$ a product structure on the complex $3$-Lie algebra $(\g_{\mathbb C},[\cdot,\cdot,\cdot]_{\g_{\mathbb C}})$. It is straightforward to see that $\g_i$ and $\g_{-i}$ are eigenspaces of $-iJ_{\mathbb C}$ corresponding to $+1$ and $-1$. Thus, $-iJ_{\mathbb C}$ is a paracomplex structure. \qed

\begin{defi}
Let $(\g,[\cdot,\cdot,\cdot]_\g)$ be a real $3$-Lie algebra. A {\bf complex product} structure on the $3$-Lie algebra $\g$ is a pair $\{J,E\}$ of a complex structure $J$ and a product structure $E$ satisfying
\begin{equation}\label{eq:compro}
 J\circ E=-E\circ J.
\end{equation}
If $E$ is perfect, we call $\{J,E\}$ a {\bf perfect complex product} structure on $\g.$
\end{defi}

\begin{rmk}
Let $\{J,E\}$ be a complex product structure on a real $3$-Lie algebra $(\g,[\cdot,\cdot,\cdot]_\g)$. For all $x\in\g_+$, by \eqref{eq:compro}, we have $E(Jx)=-Jx$, which implies that $J(\g_+)\subset\g_-$. Analogously, we obtain $J(\g_-)\subset\g_+$. Thus, we get $J(\g_-)=\g_+$ and $J(\g_+)=\g_-$. Therefore, $\dim(\g_+)=\dim(\g_-)$ and $E$ is a paracomplex structure on $\g$.
\end{rmk}

\begin{thm}
Let $(\g,[\cdot,\cdot,\cdot]_\g)$ be a real $3$-Lie algebra $(\g,[\cdot,\cdot,\cdot]_\g)$. Then the following statements are equivalent:
\begin{itemize}
\item[\rm(i)]   $\g$ has a complex product structure;
\item[\rm(ii)]  $\g$ has a complex structure $J$ and can be decomposed as $\g=\g_+\oplus\g_-$, where $\g_+,\g_-$ are $3$-Lie subalgebras of $\g$ and $\g_-=J\g_+$.
\end{itemize}
\end{thm}
\pf Let $\{J,E\}$ be a complex product structure and let $\g_{\pm}$ denote the eigenspaces corresponding to the eigenvalues $\pm1$ of $E$. By Theorem \ref{product-structure-subalgebra},  both $\g_+$ and $\g_-$
are $3$-Lie subalgebras of $\g$ and $J\circ E=-E\circ J$ implies $\g_-=J\g_+.$

Conversely, we can define a linear map $E:\g\lon\g$ by
\begin{eqnarray*}
E(x+\alpha)=x-\alpha,\,\,\,\,\forall x\in\g_+,\alpha\in\g_-.
\end{eqnarray*}
By Theorem \ref{product-structure-subalgebra}, $E$ is a product structure on $\g$. By $\g_-=J\g_+$ and $J^2=-\Id$,
we have
\begin{eqnarray*}
E(J(x+\alpha))=E(J(x)+J(\alpha))=-J(x)+J(\alpha)=-J(E(x+\alpha)).
\end{eqnarray*}
Thus, $\{J,E\}$ is a complex product structure on $\g$. The proof is finished. \qed\vspace{2mm}

\begin{ex}\label{ex:A4cp}{\rm
  Consider the product structures and the  complex structures on the $4$-dimensional Euclidean $3$-Lie algebra $A_4$ given in Example \ref{ex:A4product} and Example \ref{ex:A4complex} respectively. Then
 $\{J_i,E_i\}$ for $i=1,2,3,4,5,6$ are complex product structures on $A_4$.
  }
\end{ex}

 We give a  characterization of a   perfect complex product structure on a 3-Lie algebra.

\begin{pro}\label{3-pre-Lie-complex-product}
Let $E$ be a perfect paracomplex structure on a real $3$-Lie algebra $(\g,[\cdot,\cdot,\cdot]_\g)$. Then there is a perfect complex product structure $\{J,E\}$ on $\g$ if and only if there exists a linear isomorphism $\phi:\g_+\lon\g_-$   satisfying the following equation
\begin{eqnarray}\label{complex-perfect-product}
\nonumber\phi[x,y,z]_\g&=&-[\phi(x),\phi(y),\phi(z)]_\g+[\phi(x),y,z]_\g+[x,\phi(y),z]_\g+[x,y,\phi(z)]_\g\\
&&+\phi[\phi(x),\phi(y),z]_\g+\phi[x,\phi(y),\phi(z)]_\g+\phi[\phi(x),y,\phi(z)]_\g,\quad \forall x,y,z\in\g_+.
\end{eqnarray}
\end{pro}
\pf Let $\{J,E\}$ be a perfect complex product structure on $\g$. Define a linear isomorphism $\phi:\g_+\lon\g_-$ by $\phi\triangleq J|_{\g_+}:\g_+\lon\g_-$. By the compatibility condition \eqref{complex-structure} that the complex structure $J$ satisfies and the coherence condition \eqref{eq:coherenceconPP} that  a perfect product structure  $E$ satisfies, we deduce that
\eqref{complex-perfect-product} holds.

Conversely, we define an endomorphism $J$ of $\g$ by
\begin{eqnarray}\label{complex-perfect-product-structure}
J(x+\alpha)=-\phi^{-1}(\alpha)+\phi(x),\,\,\,\,\forall x\in\g_+,\alpha\in\g_-.
\end{eqnarray}
 It is obvious that $J$ is  an almost complex structure on $\g$ and $J\circ E=-E\circ J$. For all $\alpha,\beta,\gamma\in\g_-$, let $x,y,z\in\g_+$ such that $\phi(x)=\alpha,\phi(y)=\beta$ and $\phi(z)=\gamma$. By \eqref{complex-perfect-product} and  \eqref{eq:coherenceconPP}, we have
\begin{eqnarray*}
&&-[J\alpha,J\beta,J\gamma]_\g+[J\alpha,\beta,\gamma]_\g+[\alpha,J\beta,\gamma]_\g+[\alpha,\beta,J\gamma]_\g\\
&&+J[J\alpha,J\beta,\gamma]_\g+J[\alpha,J\beta,J\gamma]_\g+J[J\alpha,\beta,J\gamma]_\g\\
&=&[x,y,z]_\g-[x,\phi(y),\phi(z)]_\g-[\phi(x),y,\phi(z)]_\g-[\phi(x),\phi(y),z]_\g\\
&&-\phi^{-1}[x,y,\phi(z)]_\g-\phi^{-1}[\phi(x),y,z]_\g-\phi^{-1}[x,\phi(y),z]_\g\\
&=&-\phi^{-1}[\phi(x),\phi(y),\phi(z)]_\g\\
&=&J[\alpha,\beta,\gamma]_\g,
\end{eqnarray*}
which implies that \eqref{complex-structure} holds for all $\alpha,\beta,\gamma\in\g_-$. Similarly, we can deduce that \eqref{complex-structure} holds for all the other cases.
Thus,   $J$ is a complex structure and $\{J,E\}$ is a perfect complex product structure on the 3-Lie algebra $\g$.   \qed\vspace{3mm}

At the end of this section, we construct perfect complex product structure using 3-pre-Lie algebras.

\emptycomment{
\begin{pro}
Let \{J,E\} be a complex product structure on the $3$-Lie algebra $\g$ and let $(\g,\g_+,\g_-)$ be the associated double $3$-Lie algebra. Then the following assertions are equivalent:
\begin{itemize}
\item[\rm(i)]   $J$ is an abelian complex structure.
\item[\rm(ii)]  The $3$-Lie subalgebra $\g_+$ and $\g_-$ are abelian.
\item[\rm(iii)]  $E$ is an abelian product structure.
\end{itemize}
\end{pro}
\pf
}

A  nondegenerate symmetric bilinear form  $\huaB\in A^*\otimes A^*$ on  a real $3$-pre-Lie algebra $(A,\{\cdot,\cdot,\cdot\})$ is called {\bf invariant} if
\begin{eqnarray}\label{3-pre-Lie-symmetric-bilinear}
\huaB(\{x,y,z\},w)=-\huaB(z,\{x,y,w\}),\,\,\,\,\forall x,y,z,w\in A.
\end{eqnarray}
 Then $\huaB$ induces a linear isomorphism $\huaB^{\sharp}:A\lon A^*$ by
 \begin{eqnarray}
\langle\huaB^{\sharp}(x),y\rangle=\huaB(x,y),\,\,\,\,\forall x,y\in A.
\end{eqnarray}

\begin{pro}\label{pro:compro}
Let $(A,\{\cdot,\cdot,\cdot\})$ be a real $3$-pre-Lie algebra with a nondegenerate symmetric bilinear from $\huaB$.
Then there is a perfect complex product structure $\{J,E\}$ on the semidirect product $3$-Lie algebra  $ A^c\ltimes_{L^*}A^*$, where $E$ is given by \eqref{eq:defiE} and  the complex structure $J$ is given as follows:
\begin{eqnarray}\label{3-pre-Lie-complex}
J(x+\alpha)=-{\huaB^{\sharp}}^{-1}(\alpha)+\huaB^{\sharp}(x),\,\,\,\,\forall x\in A,\alpha\in A^*.
\end{eqnarray}
\emptycomment{In particular, if the bilinear form $\huaB$ is a symmetric and positive definite, the perfect complex product structure $J$ is given as follows:
\begin{eqnarray}
J(x+y^*)=-y+x^*,\,\,\,\,\forall x,y\in A,
\end{eqnarray}
where for any $x=\sum_{i=1}^{n}l_ie_i\in A$, $x^*\in A^*$ is given by  $x^*=\sum_{i=1}^{n}l_ie_i^*\in A^*.$ Here, $\{e_1,\cdots,e_n\}$ is a basis of $A$ such that $\huaB(e_i,e_j)=\delta_{ij}$ and $e_1^*,\cdots,e_n^*$ is the dual basis of $A^*$.
}
\end{pro}
\pf By Proposition \ref{paracomplex-3-pre-Lie}, $E$ is a perfect product structure on $A^c\ltimes_{L^*}A^*$. For all $x,y,z\in A^c$, we have
\begin{eqnarray*}
&&-[\huaB^{\sharp}(x),\huaB^{\sharp}(y),\huaB^{\sharp}(z)]_{L^*}+[\huaB^{\sharp}(x),y,z]_{L^*}+[x,\huaB^{\sharp}(y),z]_{L^*}+[x,y,\huaB^{\sharp}(z)]_{L^*}\\
&&+\huaB^{\sharp}[\huaB^{\sharp}(x),\huaB^{\sharp}(y),z]_{L^*}+\huaB^{\sharp}[x,\huaB^{\sharp}(y),\huaB^{\sharp}(z)]_{L^*}+\huaB^{\sharp}[\huaB^{\sharp}(x),y,\huaB^{\sharp}(z)]_{L^*}\\
&=&[\huaB^{\sharp}(x),y,z]_{L^*}+[x,\huaB^{\sharp}(y),z]_{L^*}+[x,y,\huaB^{\sharp}(z)]_{L^*}\\
&=&L^*(x,y)\huaB^{\sharp}(z)+L^*(y,z)\huaB^{\sharp}(x)+L^*(z,x)\huaB^{\sharp}(y).
\end{eqnarray*}
By \eqref{3-pre-Lie-symmetric-bilinear}, we have
\begin{eqnarray*}
\langle \huaB^{\sharp}[x,y,z]_C,w\rangle&=&\langle \huaB^{\sharp}\{x,y,z\},w\rangle+\langle \huaB^{\sharp}\{y,z,x\},w\rangle+\langle \huaB^{\sharp}\{z,x,y\},w\rangle\\
                                        &=&\huaB(\{x,y,z\},w)+\huaB(\{y,z,x\},w)+\huaB(\{z,x,y\},w)\\
                                        &=&-\huaB(z,\{x,y,w\})-\huaB(x,\{y,z,w\})-\huaB(y,\{z,x,w\})\\
                                        &=&-\langle\huaB^{\sharp}(z),\{x,y,w\}\rangle-\langle\huaB^{\sharp}(x),\{y,z,w\}\rangle
                                        -\langle\huaB^{\sharp}(y),\{z,x,w\}\rangle\\
                                        &=&\langle L^*(x,y)\huaB^{\sharp}(z),w\rangle+\langle L^*(y,z)\huaB^{\sharp}(x),w\rangle+\langle L^*(z,x)\huaB^{\sharp}(y),w\rangle,
\end{eqnarray*}
which implies that
$$
\huaB^{\sharp}[x,y,z]_C=L^*(x,y)\huaB^{\sharp}(z)+L^*(y,z)\huaB^{\sharp}(x)+L^*(z,x)\huaB^{\sharp}(y).
$$
Thus, we have
\begin{eqnarray*}
\huaB^{\sharp}[x,y,z]_C&=&-[\huaB^{\sharp}(x),\huaB^{\sharp}(y),\huaB^{\sharp}(z)]_{L^*}+[\huaB^{\sharp}(x),y,z]_{L^*}+[x,\huaB^{\sharp}(y),z]_{L^*}
+[x,y,\huaB^{\sharp}(z)]_{L^*}\\
&&+\huaB^{\sharp}[\huaB^{\sharp}(x),\huaB^{\sharp}(y),z]_{L^*}+\huaB^{\sharp}[x,\huaB^{\sharp}(y),\huaB^{\sharp}(z)]_{L^*}+\huaB^{\sharp}[\huaB^{\sharp}(x),y,\huaB^{\sharp}(z)]_{L^*}.
\end{eqnarray*}
By Proposition \ref{3-pre-Lie-complex-product}, we obtain that $\{J,E\}$ is a perfect complex product structure on   $ A^c\ltimes_{L^*}A^*.$
\emptycomment{
If the invariant bilinear form $\huaB$ is symmetric and positive definite, then $\huaB^{\sharp}(e_i)=e_i^*$, where $\{e_1,\cdots,e_n\}$ is a basis of $A$ such that $\huaB(e_i,e_j)=\delta_{ij}$ and $e_1^*,\cdots,e_n^*$ is its dual basis. The proof is finished. }\qed\vspace{3mm}

Let $(A,\{\cdot,\cdot,\cdot\})$ be a real $3$-pre-Lie algebra. On the real $3$-Lie algebra $\aff(A)=A^c\ltimes_L A$, we consider two endomorphisms $J$ and $E$ given by
\begin{eqnarray}
J(x,y)=(-y,x),\,\,\,\,E(x,y)=(x,-y),\,\,\,\,\forall x,y\in A.
\end{eqnarray}

\begin{pro}
With the above notations,
  $\{J,E\}$ is a perfect complex product structure on the $3$-Lie algebra $\aff(A)$.
\end{pro}
\pf It is obvious that $E$ is a perfect product structure on $\aff(A)$. Moreover, we have $J^2=-\Id$ and $J\circ E=-E\circ J$. Obviously $\aff(A)_+=\{(x,0)|x\in A\}, \aff(A)_-=\{(0,y)|y\in A\}$. Define $\phi:\aff(A)_+\lon\aff(A)_-$ by $\phi\triangleq J|_{\aff(A)_+}:\aff(A)_+\lon\aff(A)_-$. More precisely, $\phi(x,0)=(0,x)$. Then for all $(x,0),(y,0),(z,0)\in\aff(A)_+$, we have
\begin{eqnarray*}
&&-[\phi(x,0),\phi(y,0),\phi(z,0)]_L+[\phi(x,0),(y,0),(z,0)]_L+[(x,0),\phi(y,0),(z,0)]_L\\
&&+[(x,0),(y,0),\phi(z,0)]_L
+\phi[\phi(x,0),\phi(y,0),(z,0)]_L+\phi[(x,0),\phi(y,0),\phi(z,0)]_L\\
&&+\phi[\phi(x,0),(y,0),\phi(z,0)]_L\\
&=&[\phi(x,0),(y,0),(z,0)]_L+[(x,0),\phi(y,0),(z,0)]_L+[(x,0),(y,0),\phi(z,0)]_L\\
&=&(0,\{y,z,x\})+(0,\{z,x,y\})+(0,\{x,y,z\})\\
&=&\phi[(x,0),(y,0),(z,0)]_L.
\end{eqnarray*}
By Proposition \ref{3-pre-Lie-complex-product},   $\{J,E\}$ is a perfect complex product structure on the $3$-Lie algebra $\aff(A)$.   \qed

\section{Para-K\"{a}hler structures on $3$-Lie algebras}
In this section, we add a compatibility condition between a symplectic  structure and a paracomplex structure on a 3-Lie algebra to introduce the notion of a para-K\"{a}hler structure on a $3$-Lie algebra. A para-K\"{a}hler structure gives rise to a pseudo-Riemannian structure. We introduce the notion of a Livi-Civita product associated to a pseudo-Riemannian 3-Lie algebra and give its precise formulas using the decomposition of the original 3-Lie algebra.

\begin{defi}
Let $\omega$ be a symplectic structure and $E$ a paracomplex structure on a $3$-Lie algebra $(\g,[\cdot,\cdot,\cdot]_\g)$. The triple $(\g,\omega,E)$ is called a {\bf para-Kähler} $3$-Lie algebra if the following equality holds:
 \begin{equation}\label{eq:pk}
\omega(Ex,Ey)=-\omega(x,y),\quad \forall x,y\in\g.
\end{equation}
If $E$ is perfect, we call $(\g,\omega,E)$ a {\bf perfect para-Kähler} $3$-Lie algebra.
\end{defi}

\begin{pro}
Let $(A,\{\cdot,\cdot,\cdot\})$ be a $3$-pre-Lie algebra. Then $(A^c\ltimes_{L^*}A^*,\omega,E)$ is a perfect para-Kähler $3$-Lie algebra,  where $\omega$ is given by \eqref{phase-space} and $E$ is defined by \eqref{eq:defiE}.
\end{pro}
\pf By Theorem  \ref{3-pre-Lie-phase-space}, $(A^c\ltimes_{L^*}A^*,\omega)$ is a symplectic 3-Lie algebra. By Proposition \ref{paracomplex-3-pre-Lie}, $E$ is a perfect paracomplex structure on the phase space $T^*A^c$. For all $x_1,x_2\in A,\alpha_1,\alpha_2\in A^*$, we have
\begin{eqnarray*}
\omega(E(x_1+\alpha_1),E(x_2+\alpha_2))&=&\omega(x_1-\alpha_1,x_2-\alpha_2)=\langle -\alpha_1, x_2\rangle-\langle -\alpha_2, x_1\rangle\\
                                 &=&-\omega(x_1+\alpha_1,x_2+\alpha_2).
\end{eqnarray*}
Therefore, $(T^*A^c=A^c\ltimes_{L^*}A^*,\omega,E)$ is a perfect paraKähler 3-Lie algebra.  \qed\vspace{3mm}

 Similar as the case of   para-Kähler Lie algebras, we have the following equivalent description of a para-Kähler  $3$-Lie algebra.
\begin{thm}
Let $(\g,\omega)$ be a symplectic $3$-Lie algebra. Then there exists a paracomplex structure $E$ on the $3$-Lie algebra $(\g,[\cdot,\cdot,\cdot]_\g)$ such that $(\g,\omega,E)$ is a para-Kähler $3$-Lie algebra  if and only if there exist two isotropic $3$-Lie subalgebras $\g_+$ and $\g_-$ such that $\g=\g_+\oplus\g_-$ as the direct sum of vector spaces.
 \end{thm}
\pf Let $(\g,\omega,E)$ be a para-Kähler $3$-Lie algebra. Since $E$ is a paracomplex structure on $\g$, we have
$
\g=\g_+\oplus\g_-,
$
 where $\g_+$ and $\g_-$ are $3$-Lie subalgebras of $\g$. For all $x_1,x_2\in\g_+$, by \eqref{eq:pk}, we have
\begin{eqnarray*}
\omega(Ex_1,Ex_2)=\omega(x_1,x_2)=-\omega(x_1,x_2),
\end{eqnarray*}
which implies that $\omega(\g_+,\g_+)=0$. Thus, $\g_+$ is isotropic. Similarly, $\g_-$ is also isotropic.

Conversely, since  $\g_+$ and $\g_-$ are subalgebras, $\g=\g_+\oplus \g_-$ as vector spaces, there is a product structure $E$ on $\g$ defined by
\eqref{eq:productE}.
Moreover, since $\g=\g_+\oplus \g_-$ as vector spaces and both $\g_+$ and $\g_-$ are isotropic, we obtain that dim $\g_+$=dim $\g_-$. Thus, $E$ is a paracomplex structure on $\g$. For all $x_1,x_2\in\g_+,\alpha_1,\alpha_2\in\g_-$,  since $\g_+$ and $\g_-$ are isotropic, we have
\begin{eqnarray*}
\omega(E(x_1+\alpha_1),E(x_2+\alpha_2))&=&\omega(x_1-\alpha_1,x_2-\alpha_2)=-\omega(x_1,\alpha_2)-\omega(\alpha_1,x_2)\\
                             &=&-\omega(x_1+\alpha_1,x_2+\alpha_2).
\end{eqnarray*}
Thus,   $(\g,\omega,E)$ is a para-Kähler $3$-Lie algebra. The proof is finished. \qed

\begin{ex}\label{ex:A4pK}{\rm
  Consider the symplectic structures and the perfect paracomplex structures on the $4$-dimensional Euclidean $3$-Lie algebra $A_4$ given in Example \ref{ex:A4symplectic} and Example \ref{ex:A4product} respectively. Then
 $\{\omega_i,E_i\}$ for $i=1,2,3,4,5,6$ are perfect para-Kähler structures on $A_4$.
  }
\end{ex}

\begin{ex}\label{ex:standardpK}{\rm
Let $(\frkh,[\cdot,\cdot,\cdot]_\h$ be a 3-Lie algebra and $(\h\oplus \h^*,\omega)$ its (perfect) phase space, where $\omega$ is given by \eqref{phase-space}. Then $E:\h\oplus \h^*\longrightarrow\h\oplus \h^*$ defined by
\begin{equation}
  \label{eq:Ephasespace}
  E(x+\alpha)=x-\alpha,\quad \forall x\in\h,\alpha\in\h^*,
\end{equation}
is a (perfect) paracomplex structure and $(\h\oplus \h^*,\omega,E)$ is a (perfect) para-Kähler $3$-Lie algebra.
}
\end{ex}

 Let $(\g,\omega,E)$ be a para-Kähler $3$-Lie algebra. Then it is obvious that $\g_-$ is isomorphic to $\g_+^*$ via the symplectic structure $\omega$. Moreover, it is straightforward to deduce that
 \begin{pro}\label{pro:standardpK}
    Any para-Kähler $3$-Lie algebra   is isomorphic to the para-Kähler $3$-Lie algebra associated to a phase space of a $3$-Lie algebra.
 \end{pro}

In the sequel, we study the Levi-Civita product associated to a perfect para-Kähler 3-Lie algebra.
\begin{defi}
A {\bf pseudo-Riemannian $3$-Lie algebra} is a $3$-Lie algebra $(\g,[\cdot,\cdot,\cdot]_\g)$ endowed with a  nondegenerate symmetric bilinear form $S$. The associated Levi-Civita product is the product on $\g$, $\nabla:\otimes^3\g\longrightarrow\g$ with $(x,y,z)\longmapsto \nabla_{x,y}z$, given by the following formula:
\begin{eqnarray}\label{Levi-Civita product}
3S(\nabla_{x,y}z,w)=S([x,y,z]_\g,w)-2S([x,y,w]_\g,z)+S([y,z,w]_\g,x)+S([z,x,w]_\g,y).
\end{eqnarray}
\end{defi}

\begin{pro}
Let $(\g,S)$ be a pseudo-Riemannian $3$-Lie algebra. Then the Levi-Civita product $\{\cdot,\cdot,\cdot\}$ satisfies the following equations:
\begin{eqnarray}
\nabla_{x,y}z&=&-\nabla_{y,x}z,\\
\nabla_{x,y}z+\nabla_{y,z}x+\nabla_{z,x}y&=&[x,y,z]_\g.
\end{eqnarray}
\end{pro}
\pf For all $w\in\g$, it is obvious that
\begin{eqnarray*}
3S(\nabla_{y,x}z,w)&=&S([y,x,z]_\g,w)-2S([y,x,w]_\g,z)+S([x,z,w]_\g,y)+S([z,y,w]_\g,x)\\
               &=&-3S(\nabla_{x,y}z,w).
\end{eqnarray*}
By the nondegeneracy of $S$, we obtain $\nabla_{x,y}z=-\nabla_{y,x}z.$

For all $x,y,z,w\in\g$,   we have
\begin{eqnarray*}
3S(\nabla_{x,y}z,w)&=&S([x,y,z]_\g,w)-2S([x,y,w]_\g,z)+S([y,z,w]_\g,x)+S([z,x,w]_\g,y),\\
3S(\nabla_{y,z}x,w)&=&S([y,z,x]_\g,w)-2S([y,z,w]_\g,x)+S([z,x,w]_\g,y)+S([x,y,w]_\g,z),\\
3S(\nabla_{z,x}y,w)&=&S([z,x,y]_\g,w)-2S([z,x,w]_\g,y)+S([x,y,w]_\g,z)+S([y,z,w]_\g,x).
\end{eqnarray*}
Add up the three equations, we have
$$3S(\nabla_{x,y}z+\nabla_{y,z}x+\nabla_{z,x}y,w)=3S([x,y,z]_\g,w),$$
which implies that $\nabla_{x,y}z+\nabla_{y,z}x+\nabla_{z,x}y=[x,y,z]_\g.$ The proof is finished. \qed\vspace{3mm}

Let $(\g,\omega,E)$ be a perfect para-Kähler $3$-Lie algebra. Define a bilinear form $S$ on $\g$ by
\begin{eqnarray}
S(x,y)\triangleq \omega(x,Ey),\,\,\,\,\forall x,y\in\g.
\end{eqnarray}
\begin{pro}
With the above notations, $(\g,S)$ is a pseudo-Riemannian $3$-Lie algebra. Moreover, the associated Levi-Civita product $\nabla$ and the perfect paracomplex structure $E$ satisfy the following compatibility condition:
\begin{equation}
  E\nabla _{x,y}z=\nabla_{Ex,Ey}Ez.
\end{equation}
\end{pro}
\pf Since $\omega$ is skewsymmetric and $\omega(Ex,Ey)=-\omega(x,y)$, we have
\begin{eqnarray*}
S(y,x)=\omega(y,Ex)=-\omega(Ey,E^2x)=-\omega(Ey,x)=\omega(x,Ey)=S(x,y),
\end{eqnarray*}
which implies that $S$ is symmetric. Moreover, since $\omega$ is nondegenerate and $E^2=\Id$, it is obvious that $S$ is nondegenerate. Thus, $S$ is a pseudo-Riemannian metric on the 3-Lie algebra $\g$.  Moreover, we have
\begin{eqnarray*}
&&3S(\nabla_{Ex,Ey}Ez,w)\\&=&S([Ex,Ey,Ez]_\g,w)-2S([Ex,Ey,w]_\g,Ez)+S([Ey,Ez,w]_\g,Ex)+S([Ez,Ex,w]_\g,Ey)\\
&=& S(E[x,y,z]_\g,w)-2S(E[x,y,Ew]_\g,Ez)+S(E[y,z,Ew]_\g,Ex)+S(E[z,x,Ew]_\g,Ey)\\
&=&-( S([x,y,z]_\g,Ew)-2S([x,y,Ew]_\g,z)+S([y,z,Ew]_\g,x)+S([z,x,Ew]_\g,y))\\
&=&-3S(\nabla_{x,y}z,Ew)\\
&=&3S(E\nabla_{x,y}z,w).
\end{eqnarray*}
Thus, we have  $ E\nabla _{x,y}z=\nabla_{Ex,Ey}Ez.$\qed\vspace{3mm}

The following two propositions clarifies the relationship between the Levi-Civita product and the 3-pre-Lie multiplication on a para-Kähler $3$-Lie algebra.

\begin{pro}\label{Levi-Civita-3-pre-Lie}
Let $(\g,\omega,E)$ be a para-Kähler $3$-Lie algebra and $\nabla$  the associated Levi-Civita product.    Then for all $x_1,x_2,x_3\in\g_+$ and $\alpha_1,\alpha_2,\alpha_3\in\g_-$, we have
$$
\nabla_{x_1,x_2}x_3=\{x_1,x_2,x_3\},\quad \nabla_{\alpha_1,\alpha_2}\alpha_3=\{\alpha_1,\alpha_2,\alpha_3\}.
$$
\end{pro}
\pf Since $(\g,\omega,E)$ is a para-Kähler $3$-Lie algebra, $3$-Lie subalgebras $\g_+$ and $\g_-$ are isotropic and $\g=\g_+\oplus\g_-$ as vector spaces. For all $x_1,x_2,x_3,x_4\in\g_+$, we have
\begin{eqnarray*}
&&3\omega(\nabla_{x_1,x_2}x_3,x_4)\\&=&3S(\nabla_{x_1,x_2}x_3,Ex_4)=3S(\nabla_{x_1,x_2}x_3,x_4)\\
                            &=&S([x_1,x_2,x_3]_\g,x_4)-2S([x_1,x_2,x_4]_\g,x_3)+S([x_2,x_3,x_4]_\g,x_1)+S([x_3,x_1,x_4]_\g,x_2)\\
                            &=&\omega([x_1,x_2,x_3]_\g,x_4)-2\omega([x_1,x_2,x_4]_\g,x_3)+\omega([x_2,x_3,x_4]_\g,x_1)+\omega([x_3,x_1,x_4]_\g,x_2)\\
                            &=&0.
\end{eqnarray*}
By $(\g_+)^{\perp}=\g_+$, we obtain $\nabla_{x_1,x_2}x_3\in\g_+$. Similarly, for all $\alpha_1,\alpha_2,\alpha_3\in\g_-$,  $\nabla_{\alpha_1,\alpha_2}\alpha_3\in\g_-$.
Furthermore, for all $x_1,x_2,x_3\in\g_+,$ and $\alpha\in\g_-$, we have
\begin{eqnarray*}
&&3\omega(\nabla_{x_1,x_2}x_3,\alpha)\\&=&3S(\nabla_{x_1,x_2}x_3,E\alpha)=-3S(\nabla_{x_1,x_2}x_3,\alpha)\\
                            &=&-S([x_1,x_2,x_3]_\g,\alpha)+2S([x_1,x_2,\alpha]_\g,x_3)-S([x_2,x_3,\alpha]_\g,x_1)-S([x_3,x_1,\alpha]_\g,x_2)\\
                            &=&\omega([x_1,x_2,x_3]_\g,\alpha)+2\omega([x_1,x_2,\alpha]_\g,x_3)-\omega([x_2,x_3,\alpha]_\g,x_1)-\omega([x_3,x_1,\alpha]_\g,x_2)\\
                            &=&\omega([\alpha,x_1,x_2]_\g,x_3)+2\omega([x_1,x_2,\alpha]_\g,x_3)\\
                            &=&-3\omega(x_3,[x_1,x_2,\alpha]_\g)\\
                            &=&3\omega(\{x_1,x_2,x_3\},\alpha).
\end{eqnarray*}
Thus,  $\nabla_{x_1,x_2}x_3=\{x_1,x_2,x_3\}.$ Similarly, we have  $\nabla_{\alpha_1,\alpha_2}\alpha_3=\{\alpha_1,\alpha_2,\alpha_3\}$. The proof is finished.\qed

\begin{pro}
Let $(\g,\omega,E)$ be a perfect para-Kähler $3$-Lie algebra and $\nabla$  the associated Levi-Civita product. Then for all $x_1,x_2\in\g_+$ and $\alpha_1,\alpha_2\in\g_-$, we have
\begin{eqnarray}
\label{eq:conn1}\nabla_{x_1,x_2}\alpha_1&=&\{x_1,x_2,\alpha_1\}+\frac{2}{3}(\{x_2,\alpha_1,x_1\}+\{\alpha_1,x_1,x_2\}),\\
\label{eq:conn2}\nabla_{\alpha_1,x_1}x_2&=&-\frac{1}{3}\{\alpha_1,x_1,x_2\}+\frac{2}{3}\{x_2,\alpha_1,x_1\},\\
\label{eq:conn3}\nabla_{\alpha_1,\alpha_2}x_1&=&\{\alpha_1,\alpha_2,x_1\}+\frac{2}{3}(\{\alpha_2,x_1,\alpha_1\}+\{x_1,\alpha_1,\alpha_2\}),\\
\label{eq:conn4}\nabla_{x_1,\alpha_1}\alpha_2&=&-\frac{1}{3}\{x_1,\alpha_1,\alpha_2\}+\frac{2}{3}\{\alpha_2,x_1,\alpha_1\}.
\end{eqnarray}
\end{pro}
\pf Since $(\g,\omega,E)$ is a perfect para-Kähler $3$-Lie algebra, $3$-Lie subalgebras $\g_+$ and $\g_-$ are isotropic and $\g=\g_+\oplus\g_-$ as vector spaces. Thus, we have $S(\g_+,\g_+)=S(\g_-,\g_-)=0.$ For all $x_1,x_2\in\g_+$ and $\alpha_1,\alpha_2\in\g_-$, we have
\begin{eqnarray*}
&&3S(\nabla_{x_1,x_2}\alpha_1,\alpha_2)\\
&=&S([x_1,x_2,\alpha_1]_\g,\alpha_2)-2S([x_1,x_2,\alpha_2]_\g,\alpha_1)
+S([x_2,\alpha_1,\alpha_2]_\g,x_1)+S([\alpha_1,x_1,\alpha_2]_\g,x_2)=0.
\end{eqnarray*}
 Since $S$ is nondegenerate, we have $\nabla_{x_1,x_2}\alpha_1\in\g_-.$ Moreover, For all $x_1,x_2,x_3\in\g_+$ and $\alpha_1\in\g_-$, we have
\begin{eqnarray*}
&&3\omega(\nabla_{x_1,x_2}\alpha_1,x_3)\\&=&3S(\nabla_{x_1,x_2}\alpha_1,Ex_3)=3S(\nabla_{x_1,x_2}\alpha_1,x_3)\\
                                     &=&S([x_1,x_2,\alpha_1]_\g,x_3)-2S([x_1,x_2,x_3]_\g,\alpha_1)
+S([x_2,\alpha_1,x_3]_\g,x_1)+S([\alpha_1,x_1,x_3]_\g,x_2)\\
                                     &=&\omega([x_1,x_2,\alpha_1]_\g,Ex_3)-2\omega([x_1,x_2,x_3]_\g,E\alpha_1)
+\omega([x_2,\alpha_1,x_3]_\g,Ex_1)+\omega([\alpha_1,x_1,x_3]_\g,Ex_2)\\
                                     &=&\omega([x_1,x_2,\alpha_1]_\g,x_3)+2\omega([x_1,x_2,x_3]_\g,\alpha_1)
+\omega([x_2,\alpha_1,x_3]_\g,x_1)+\omega([\alpha_1,x_1,x_3]_\g,x_2)\\               &=&\omega([x_1,x_2,\alpha_1]_\g,x_3)+2\omega(\{x_1,x_2,\alpha_1\},x_3)
                                     +\omega(\{x_2,\alpha_1,x_1\},x_3)+\omega(\{\alpha_1,x_1,x_2\},x_3).
\end{eqnarray*}
Thus, we obtain
\begin{eqnarray*}
\nabla_{x_1,x_2}\alpha_1&=&\{x_1,x_2,\alpha_1\}+\frac{2}{3}(\{x_2,\alpha_1,x_1\}+\{\alpha_1,x_1,x_2\}),
\end{eqnarray*}
which implies that \eqref{eq:conn1} holds.

For all $x_1,x_2\in\g_+$ and $\alpha_1,\alpha_2\in\g_-$, we have
\begin{eqnarray*}
&&3S(\nabla_{\alpha_1,x_1}x_2,\alpha_2)\\
&=&S([\alpha_1,x_1,x_2]_\g,\alpha_2)-2S([\alpha_1,x_1,\alpha_2]_\g,x_2)
+S([x_1,x_2,\alpha_2]_\g,\alpha_1)+S([x_2,\alpha_1,\alpha_2]_\g,x_1)=0.
\end{eqnarray*}
Since $S$ is nondegenerate, we have $\nabla_{\alpha_1,x_1}x_2\in\g_-.$ Moreover, For all $x_1,x_2,x_3\in\g_+$ and $\alpha_1\in\g_-$, we have
\begin{eqnarray*}
&&3\omega(\nabla_{\alpha_1,x_1}x_2,x_3)\\&=&3S(\nabla_{\alpha_1,x_1}x_2,Ex_3)=3S(\nabla_{\alpha_1,x_1}x_2,x_3)\\
                                     &=&S([\alpha_1,x_1,x_2]_\g,x_3)-2S([\alpha_1,x_1,x_3]_\g,x_2)
+S([x_1,x_2,x_3]_\g,\alpha_1)+S([x_2,\alpha_1,x_3]_\g,x_1)\\                                     &=&\omega([\alpha_1,x_1,x_2]_\g,x_3)-2\omega([\alpha_1,x_1,x_3]_\g,x_2)
-\omega([x_1,x_2,x_3]_\g,\alpha_1)+\omega([x_2,\alpha_1,x_3]_\g,x_1)\\
 &=&\omega([\alpha_1,x_1,x_2]_\g,x_3)-2\omega(\{\alpha_1,x_1,x_2\},x_3)
-\omega(\{x_1,x_2,\alpha_1\},x_3)+\omega(\{x_2,\alpha_1,x_1\},x_3).
\end{eqnarray*}
Thus, we obtain
\begin{eqnarray*}
\nabla_{\alpha_1,x_1}x_2&=&-\frac{1}{3}\{\alpha_1,x_1,x_2\}+\frac{2}{3}\{x_2,\alpha_1,x_1\},
\end{eqnarray*}
which implies that \eqref{eq:conn2} holds.

  \eqref{eq:conn3}  and   \eqref{eq:conn4} can be proved similarly. We omit details. The proof is finished. \qed\vspace{3mm}

Under the isomorphism given in Proposition \ref{pro:standardpK} and the correspondence given in Theorem \ref{thm:MT-ps},  using the formulas provided in Proposition \ref{pro:stuctureMP3preLie}, we get
\begin{cor}
 For the perfect para-Kähler $3$-Lie algebra $(\h\oplus \h^*,\omega,E)$ given in Example \ref{ex:standardpK},   for all $x_1,x_2\in\h$ and $\alpha_1,\alpha_2\in\h^*$, we have
\begin{eqnarray}
\label{eq:conn11}\nabla_{x_1,x_2}\alpha_1&=&(L^*(x_1,x_2)-\frac{1}{3}R^*(x_2,x_1)+\frac{1}{3}R^*(x_1,x_2))\alpha_1,\\
\label{eq:conn22}\nabla_{\alpha_1,x_1}x_2&=&(\frac{1}{3}R^*(x_1,x_2)+\frac{2}{3}R^*(x_2,x_1))\alpha_1,\\
\label{eq:conn33}\nabla_{\alpha_1,\alpha_2}x_1&=&(\huaL^*(\alpha_1,\alpha_2)-\frac{1}{3}\huaR^*(\alpha_2,\alpha_1)+\frac{1}{3}\huaR^*(\alpha_1,\alpha_2))x_1,\\
\label{eq:conn44}\nabla_{x_1,\alpha_1}\alpha_2&=&(\frac{1}{3}\huaR^*(\alpha_1,\alpha_2)+\frac{2}{3}\huaR^*(\alpha_2,\alpha_1))x_1.
\end{eqnarray}
\end{cor}

\section{Pseudo-K\"{a}hler structures on $3$-Lie algebras}

In this section, we add a compatibility condition between a symplectic  structure and a  complex structure on a 3-Lie algebra to introduce the notion of a pseudo-K\"{a}hler structure on a $3$-Lie algebra. The relation between para-K\"{a}hler structures and    pseudo-K\"{a}hler structures on a 3-Lie algebra is investigated.
\begin{defi}
Let $\omega$ be a symplectic structure and $J$ a complex structure on a real $3$-Lie algebra $(\g,[\cdot,\cdot,\cdot]_\g)$.  The triple $(\g,\omega,J)$ is called a real {\bf pseudo-Kähler} $3$-Lie algebra if
\begin{equation}\label{eq:pK}
\omega(Jx,Jy)=\omega(x,y),\quad \forall x,y\in\g.
\end{equation}
\end{defi}

\begin{ex}\label{ex:A4sK}{\rm
  Consider the symplectic structures and the  complex structures on the $4$-dimensional Euclidean $3$-Lie algebra $A_4$ given in Example \ref{ex:A4symplectic} and Example \ref{ex:A4complex} respectively. Then
 $\{\omega_i,J_i\}$ for $i=1,2,3,4,5,6$ are pseudo-Kähler structures on $A_4$.
  }
\end{ex}

\begin{pro}
Let $(\g,\omega,J)$ be a real pseudo-Kähler $3$-Lie algebra. Define a bilinear form $S$ on $\g$ by
\begin{eqnarray}
S(x,y)\triangleq \omega(x,Jy),\,\,\,\,\forall x,y\in\g.
\end{eqnarray}
Then $(\g,S)$ is a pseudo-Riemannian $3$-Lie algebra.
\end{pro}
\pf By \eqref{eq:pK}, we have
\begin{eqnarray*}
S(y,x)=\omega(y,Jx)=\omega(Jy,J^2x)=-\omega(Jy,x)=\omega(x,Jy)=S(x,y),
\end{eqnarray*}
which implies that $S$ is symmetric. Moreover, since $\omega$ is nondegenerate and $J^2=-\Id$, it is obvious that $S$ is nondegenerate. Thus, $S$ is a pseudo-Riemannian metric on the 3-Lie algebra $\g$.   \qed

\begin{defi}
Let $(\g,\omega,J)$ be a real pseudo-Kähler $3$-Lie algebra. If the associated pseudo-Riemannian metric is positive definite, we call $(\g,\omega,J)$ a real {\bf Kähler} $3$-Lie algebra.
\end{defi}

\begin{thm}
Let $(\g,\omega,E)$ be a complex para-Kähler $3$-Lie algebra. Then $(\g_{\mathbb R},\omega_{\mathbb R},J)$ is a real pseudo-Kähler $3$-Lie algebra, where $\g_{\mathbb R}$ is the underlying
real $3$-Lie algebra, $J=iE$ and  $\omega_{\mathbb R}=\re(\omega)$ is the real part of $\omega.$
\end{thm}
\pf By Proposition \ref{equivalent}, $J=iE$ is a complex structure on the complex $3$-Lie algebra $\g$. Thus, $J$ is also a complex structure on the real $3$-Lie algebra $\g_{\mathbb R}$. It is obvious that $\omega_{\mathbb R}$ is skew-symmetric. If for all $x\in\g$,   $\omega_{\mathbb R}(x,y)=0$. Then we have
\begin{eqnarray*}
\omega(x,y)=\omega_{\mathbb R}(x,y)+i\omega_{\mathbb R}(-ix,y)=0.
\end{eqnarray*}
By the nondegeneracy of $\omega$, we obtain $y=0$. Thus, $\omega_{\mathbb R}$ is nondegenerate. Therefore, $\omega_{\mathbb R}$ is a symplectic structure on the real $3$-Lie algebra $\g_{\mathbb R}$. By $\omega(Ex,Ey)=-\omega(x,y)$, we have
\begin{eqnarray*}
\omega_{\mathbb R}(Jx,Jy)=\re(\omega(iEx,iEy))=\re(-\omega(Ex,Ey))=\re(\omega(x,y))=\omega_{\mathbb R}(x,y).
\end{eqnarray*}
Thus, $(\g_{\mathbb R},iE,\omega_{\mathbb R})$ is a real pseudo-Kähler $3$-Lie algebra.  \qed

Conversely, we have
\begin{thm}
Let $(\g,\omega,J)$ be a real pseudo-Kähler $3$-Lie algebra. Then $(\g_{\mathbb C},\omega_{\mathbb C},E)$ is a complex para-Kähler $3$-Lie algebra, where $\g_{\mathbb C}=\g\otimes_{\mathbb R}\mathbb C$ is the complexification of $\g$,  $E=-iJ_{\mathbb C}$ and $\omega_{\mathbb C}$ is the complexification of $\omega$, more precisely,
\begin{eqnarray}\label{complex-omega}
\omega_{\mathbb C}(x_1+iy_1,x_2+iy_2)=\omega(x_1,x_2)-\omega(y_1,y_2)+i\omega(x_1,y_2)+i\omega(y_1,x_2), \quad\forall x_1,x_2,y_1,y_2\in\g.
\end{eqnarray}
\end{thm}
\pf By Corollary \ref{complex-to-special-paracomplex}, $E=-iJ_{\mathbb C}$ is a paracomplex structure on the complex $3$-Lie algebra $\g_{\mathbb C}$.
It is obvious that $\omega_{\mathbb C}$ is skew-symmetric and nondegenerate. Moreover,
since $\omega$  is a symplectic structure on $\g$, we deduce that $\omega_{\mathbb C}$ is a symplectic structure on $\g_{\mathbb C}.$ Finally, by $\omega(Jx,Jy)=\omega(x,y)$, we have
\begin{eqnarray*}
\omega_{\mathbb C}(E(x_1+iy_1),E(x_2+iy_2))&=&\omega_{\mathbb C}(Jy_1-iJx_1,Jy_2-iJx_2)\\
                                           &=&\omega(Jy_1,Jy_2)-\omega(Jx_1,Jx_2)-i\omega(Jx_1,Jy_2)-i\omega(Jy_1,Jx_2)\\
                                           &=&\omega(y_1,y_2)-\omega(x_1,x_2)-i\omega(x_1,y_2)-i\omega(y_1,x_2)\\
                                           &=&-\omega_{\mathbb C}(x_1+iy_1,x_2+iy_2).
\end{eqnarray*}
Therefore, $(\g_{\mathbb C},\omega_{\mathbb C},-iJ_{\mathbb C})$ is a complex  para-Kähler $3$-Lie algebra.  \qed\vspace{3mm}

At the end of this section, we construct a Kähler $3$-Lie algebra using a $3$-pre-Lie algebra with a symmetric and positive definite invariant bilinear form.

\begin{pro}
Let $(A,\{\cdot,\cdot,\cdot\})$ be a real $3$-pre-Lie algebra with a symmetric and positive definite invariant bilinear form $\huaB$. Then $(A^c\ltimes_{L^*}A^*,\omega,-J)$ is a real Kähler $3$-Lie algebra, where $J$ is given by   \eqref{3-pre-Lie-complex} and $\omega$ is given by   \eqref{phase-space}.
\end{pro}
\pf By Theorem \ref{3-pre-Lie-phase-space} and Proposition \ref{pro:compro},   $\omega$ is a symplectic structure and $J$ is a perfect complex structure on the semidirect product 3-Lie algebra $( A^c\ltimes_{L^*}A^*,[\cdot,\cdot,\cdot]_{L^*})$. Obviously, $-J$ is also a perfect complex structure on $A^c\ltimes_{L^*}A^*$. Let
$\{e_1,\cdots,e_n\}$ be a basis of $A$ such that $\huaB(e_i,e_j)=\delta_{ij}$ and $e_1^*,\cdots,e_n^*$ be the dual basis of $A^*$. Then for all $i,j,k,l$, we have
\begin{eqnarray*}
\omega(e_i+e_j^*,e_k+e_l^*)&=&\delta_{jk}-\delta_{li},\\
\omega(-J(e_i+e_j^*),-J(e_k+e_l^*))&=&\omega(e_j-e_i^*,e_l-e_k^*)=-\delta_{il}+\delta_{kj},
\end{eqnarray*}
which implies that $\omega(-J(x+\alpha),-J(y+\beta))=\omega(x+\alpha,y+\beta)$ for all $x,y\in A$ and $\alpha,\beta\in A^*$. Therefore, $(A^c\ltimes_{L^*}A^*,\omega,-J)$ is a pseudo-Kähler $3$-Lie algebra. Finally, Let $x=\sum_{i=1}^{n}\lambda_ie_i\in A,\alpha=\sum_{i=1}^{n}\mu_ie_i^*\in A^*$ such that $x+\alpha\not=0.$   We have
\begin{eqnarray*}
S(x+\alpha,x+\alpha)&=&\omega(x+\alpha,-J(x+\alpha))\\
                    &=&\omega\big(\sum_{i=1}^{n}\lambda_ie_i+\sum_{i=1}^{n}\mu_ie_i^*,\sum_{i=1}^{n}\mu_ie_i-\sum_{i=1}^{n}\lambda_ie_i^*)\big)\\
                    &=&\sum_{i=1}^{n}\mu_i^2+\sum_{i=1}^{n}\lambda_i^2>0.
\end{eqnarray*}
Thus, $S$ is positive definite. Therefore, $\{A^c\ltimes_{L^*}A^*,\omega,-J\}$ is a real  Kähler $3$-Lie algebra. \qed

 Department of Mathematics, Jilin University,
 Changchun 130012,  China

 Email:
shengyh@jlu.edu.cn, tangrong16@mails.jlu.edu.cn

\end{document}